\documentclass[preprint,9pt]{elsarticle}

 \usepackage{graphics}

\usepackage{amssymb}

 \usepackage{amsthm}

\begin{document}

\begin{frontmatter}



\title{Efficient determination of critical parameters of nonlinear Schr\"{o}dinger equation with point-like potential using  generalized  polynomial chaos methods}
\author{Debananda Chakraborty, Jae-Hun Jung}
\address{Department of Mathematics, State University of New York at Buffalo, Buffalo, NY 14260-2900, USA.}
\ead{dc58@buffalo.edu, jaehun@buffalo.edu}
\author{Emmanuel Lorin}
\address{School of Mathematics and Statistics, Carleton University, Ottawa, Canada K1S 5B6}
\ead{elorin@math.carleton.ca}
\cortext[cor1]{Corresponding author}


\begin{abstract}
We consider the nonlinear Schr\"{o}dinger equation with a point-like source term. The soliton interaction with such a singular potential yields a critical solution behavior. That is, for the given value of the potential strength and the soliton amplitude, there exists a critical velocity of the initial soliton solution, around which the solution is 
either trapped by or transmitted through the potential. In this paper, we propose an efficient method for finding such a critical velocity by using the generalized polynomial chaos method. For the proposed method, we assume that the soliton velocity is a random variable and expand the solution in the random space using the orthogonal polynomials. The 
proposed method finds the critical velocity accurately with spectral convergence. Thus the computational complexity is much reduced. Numerical results for the smaller and higher values of the potential strength confirm the spectral convergence of the proposed method.

\end{abstract}

\begin{keyword}
Nonlinear Schr\"{o}dinger equation, Singular potential, 
Generalized polynomial chaos, Stochastic collocation method, Split step 
Fourier method, Spectral convergence
\end{keyword}

\end{frontmatter}
\section{Introduction}
Nonlinear Schr\"{o}dinger equation (NLSE) describes a broad range of physical phenomena,  e.g. nonlinear modulation of collisionless plasma waves \cite{Ablowitz}, self trapping of a light beam in a color dispersive system \cite{Agarwal}, helical motion in a very thin vortex filament \cite{Hasegawa}, propagation of heat pulses in an-harmonic crystals \cite{Drazin}, modulation instability in water waves \cite{Hasegawa}, etc.  In optical fibers, the soliton solutions of the NLSE provide a secure means to carry bits of information over many thousands of miles \cite{Agarwal}. Termed as the Gross-Pitaveskii equation, the NLSE with an appropriate potential can be utilized to describe the dynamics of the Bose-Einstein condensate, both with the attractive and repulsive nonlinearities \cite{Hasegawa, Malomed}. It is our objective in this paper to solve the Gross-Pitaveskii equation equipped with a point-like potential to find the critical values of the  soliton velocities when the amplitude of the  point-like potential is either very small $(\sim 10^{-1} - 10^{-2})$ or large $(\sim 2.5 - 4.5)$  compared to the soliton amplitude which is the unity in our paper.

We know the soliton solution of the homogeneous NLSE 
\begin{eqnarray}
i\partial_{t}u + \frac{1}{2}\partial^{2}_{x}u + u|u|^{2}  = 0, \quad -\infty < x< \infty, \; t > 0 ,
\label{NLSE}
\end{eqnarray}
with initial condition $u_0$ given by
\begin{equation}
u_{0}(x) = A \mathrm{sech}\left(A\left(x\right)\right)e^{\left(i\phi +iVx\right)}
\label{initial0}
\end{equation}
is given by
\begin{eqnarray}
u(x,t) = A \mathrm{sech} \left(A\left(x-Vt\right)\right)\exp\left(i\phi + iVx + \frac{i}{2} \left(A^{2}-V^{2}\right)t\right), \; A > 0, \; V \in \mathbb{R},
\label{NLSsol}
\end{eqnarray}
where $A$ is the soliton amplitude, $V$  the soliton velocity, and $\phi$ the phase lag. 
Consider a perturbed NLSE, that is, the Gross-Pitaveskii equation by adding an external potential, $-\epsilon \delta(x)u$,
\begin{eqnarray}
\left\{ \begin{array}{cc} i\partial_{t}u + \frac{1}{2}\partial^{2}_{x}u + u|u|^{2}  = -\epsilon \delta(x)u, \\
u(x,0) = u_{0}(x), \end{array} \right.
\label{GP1}
\end{eqnarray}
where $\delta(x)$ is the Dirac delta function with a constant $\epsilon \in \mathbb{R}$. Such an external potential represents the impurity or defect in the optical fiber. The well-posedness of the equation
\begin{eqnarray}
\left\{ \begin{array}{cc} i\partial_{t}u + \frac{1}{2}\partial^{2}_{x}u + u|u|^{p-1}  = -\epsilon \delta(x)u, \\
u(x,0) = u_{0}(x), \end{array} \right.
\label{GP2}
\end{eqnarray}
with $p\geq 1$ and initial data $u_0$ in $H^1({\ensuremath{\mathbb{R}}})$, has been extensively studied and is based on the knowledge of the self-adjoint (in $L^2$) operator $-\partial_{xx}+\epsilon\delta$. Using \cite{Cazenave}, Le Coz {\it et al.} proved the existence of a time $T>0$ and of a unique solution to Eq. \ref{GP1} (where $\epsilon \in {\ensuremath{\mathbb{R}}}$) in $C\big([0,T),H^1({\ensuremath{\mathbb{R}}})\big)\cap C^1\big([0,T),H^{-1}({\ensuremath{\mathbb{R}}})\big)$ satisfying $\lim_{t \rightarrow T} \|\partial_xu\|_2=\infty$. Moreover the energy is conserved in time. This result was extended to $p\geq 1$ by Fukuizumi {\it et al.} in \cite{Fuku2}.  For $p=3$ (more generally $p\in (1,5)$), global existence in $H^1$ also holds by Gagliardo-Nirenberg's inequality and energy conservation. Global existence in $H^1$, is also discussed by Goodman {\it et al.} \cite{GHW} using a fixed-point argument and time-invariance of the $L^2$-norm and of the Hamiltonian derived from the NLSE.  Notice that the study of stability of nonlinear bound states which are solutions of the form $\exp(-{\tt i}\omega t)\phi_{\omega}(x)$ with $\omega>0$, and for which:
\begin{eqnarray}
\label{NLSB}
-\frac{1}{2}\partial_{xx}\phi_{\omega}  -\epsilon \phi_{\omega} - |\phi_{\omega}|^2 \phi_{\omega}=\omega \phi_{\omega}
\end{eqnarray}
plays an important role in the theory of NLSE with defect and could possibly be useful numerically. Explicit formulas and stability analysis for $\phi_{\omega}$ can be found in \cite{Fuku2,LeCoz}.

If we now take a soliton approaching the impurity from the left as an initial condition $u_{0}$:
\begin{equation}
u_{0}(x) = A \mathrm{sech}\left(A\left(x-x_{0}\right)\right)e^{\left(i\phi +iVx\right)},\; x_{0} \ll 0,
\label{initial1}
\end{equation}
then until the time $t_{0} = \frac{x_{0}}{V}$, the solution will still be given by 
Eq. \ref{NLSsol}. In this paper we consider $A = 1$ and $\phi = 0$. Thus the soliton velocity $V$ and the strength of the impurity  $\epsilon$ are the only parameters of the problem.

For $t_{0} > \frac{x_{0}}{V}$, the effects of the potential are highly visible and a lot of research has been done on the transmission and reflection coefficients of the $\delta$-potential by the standard scattering theory \cite{Holmer}.
Malomed and his co-workers \cite{Cao_Malomed, Malomed} showed mainly numerically, that for any given velocity $V \left( > 0 \right)$, there exists a threshold value $\epsilon_{thr} \left( > 0 \right)$ of $\epsilon$, for which the soliton can marginally pass through the defect. So for the given velocity $V$, if $\epsilon < \epsilon_{thr}$, the soliton can pass through the defect and  the soliton gets trapped otherwise. They considered the soliton-soliton collisions within the coupled NLSE. In the limiting condition one soliton has very large amplitude and  is very narrow accordingly, while the soliton governed by the other equation has finite amplitude and width. In this limiting condition the two coupled NLSE are reduced to a single equation, in which the narrow soliton will be represented by the $\delta$-function,
 \begin{equation}
 i\partial_{t}u + \frac{1}{2}\partial^{2}_{x}u + u|u|^{2}  = -\epsilon \delta(x)u. \nonumber
 \label{NLSE2}
 \end{equation}

Holmer and his co-workers studied the NLSE with $V \gg 1$ \cite{Holmer1} and $V \ll 1,\; \epsilon \ll 1$ \cite{Holmer2}.  They showed for high $V$, there exits the bound state which is given by 
$$
u(x,t) = e^{i\lambda^{2}\frac{t}{2}}\lambda \mathrm{sech}\left(\lambda |x| + \tanh^{-1}\left(\epsilon/\lambda\right)\right), \quad 0<\lambda <\epsilon, 
$$
and this bound state is ``left behind" after the interaction (see bottom right figure of Figure \ref{fig:soliton3}).  Also they proved in \cite{Holmer2} that for $V \ll 1$ and $ \epsilon \ll 1$, the solution can be approximated by the soliton solution of the homogeneous NLSE $\left(\epsilon = 0\right)$. 
To solve Eq. \ref{NLSE2} for any given $\epsilon \left( > 0 \right)$ and $V \left( > 0 \right)$, we consider three cases: (a) small value of $\epsilon, \; \mathrm{where} \; \epsilon \le 0.3$ (b) moderate value of  $\epsilon, \; \mathrm{where} \; 0.3 < \epsilon \le 3.5$ \cite{Cao_Malomed} and (c) large value of $\epsilon, \; \mathrm{where} \; \epsilon > 3.5$. For solving Eq. \ref{NLSE2}, one can use the Split Step Fourier Method (SSFM) to reduce the computational time. To get $\epsilon_{thr}$ for any given $V$  with certain accuracy one must conduct a series of simulations. The number of simulations increase with the increase of the level of accuracy. In addition to conduct a series of simulations with small time steps, one needs  a large amount of the computational time. This is our main motivation to propose a suitable method to overcome such a high computational complexity by using the generalized polynomial chaos (gPC) methods \cite{XiuBook}.

The gPC method belongs to the class of non-sampling methods \cite{Xiu2, Xiu3}. In this method the stochastic quantities are expanded by orthogonal polynomials. Different types of orthogonal polynomials can be chosen for better convergence. The  gPC expansion is a  spectral representation in random space and exhibits fast convergence when the expanded function depends smoothly on the random parameters \cite{Gottlieb}. When the gPC method is applied to solve any differential equation, the main computational work is needed to solve the expansion coefficients of the gPC expansion. A common approach is the Galerkin method that minimizes the residue in the polynomial space. The stochastic Galerkin (SG) approach, however, would be extremely difficult to use when the governing stochastic equations take complicated forms. In our case, the NLSE contains the nonlinear term $|u|^{2}u$. For the SG method, it is very hard to get the corresponding explicit deterministic equations after expanding the nonlinear terms. So that, in this work we use the high-order stochastic collocation (SC) approach \cite{Xiu2} that combines the advantages of both the Monte Carlo sampling and the gPC-Galerkin methods. The gPC method reduces the number of simulations for finding the critical velocity, $V_{c}$, for any given value of $\epsilon$ thanks to the high-order convergence of the method. Since the equation has only two parameters, i.e. $\epsilon$ and $V$, we treat at least one of them as a stochastic variable in the gPC framework. In the present work we consider $V$ as the stochastic variable and let $\epsilon$ be fixed.  So for any given $\epsilon$, we find $V_{c}$, the critical value of $V$ around which the soliton is either transmitted or trapped. Thus it is obvious that for   $V > V_{c}$, the soliton passes through the defect. By adopting this idea we develop a step-by-step gPC collocation method to find the critical velocity of the soliton.

In  \cite{Cao_Malomed} the relation between $\epsilon_{thr}$ and $V$ was obtained only for the moderate values of $\epsilon$, i.e. for those  comparable to the soliton amplitude $A = 1$. But the results of the numerical simulations for very small  or large values of $V$ were not obtained, perhaps due to the huge computational burden. By the gPC method, we were able to reduce the overhead computational time, for having detailed simulations  performed for large and small values of $\epsilon$ to find the corresponding critical velocity $V_{c}$.

Since the analysis for the moderate values of $\epsilon$ are already done \cite{Cao_Malomed}, we do not intended to repeat the analysis for those values of $\epsilon$ in this paper. Here we mainly focus on the small and high values of $\epsilon$. For the small values of $\epsilon$, the gPC takes much longer time than the gPC method for the large values of $\epsilon$ due to the extremely small critical velocities. 

This paper is organized as follows. In Section $2$, we discuss the SSFM. Section $3$ describes the  gPC collocation method. Section $4$ contains the gPC collocation algorithm for the NLSE with the singular potential term to detect the critical velocity for the given value of $\epsilon$. Section $5$ presents the numerical results. Concluding remarks and future works are presented in Section $6$.

\section{Split Step Fourier Method}
The SSFM is a pseudo-spectral numerical method used to solve nonlinear PDEs like the NLSE. Eq. \ref{NLSE} can be rewritten as 
\begin{eqnarray}
\frac{\partial u}{\partial t} = i\left[N + D\right]u,
\label{split1}
\end{eqnarray}
where $D = \frac{1}{2}\frac{\partial^{2}}{\partial x^{2}}$ and $N = |u|^{2}$.
The solution of Eq. \ref{split1} can be written as
$$
u(x,t) = e^{it\left(D+N\right)}u(x,0),
$$
where $u(x,0)$ is the initial condition. Since $D$ and $N$ are the operators, they do not necessarily commute. However the Baker-Hausdorff formula can be applied to show that the error  will be of order $dt^{2}$ if we are taking a small but finite time step $dt$ \cite{sinkin}.  We therefore can write 
\begin{equation}
u\left(x,t+dt\right) \approx e^{idtN}e^{idtD}u(x,t).
\label{split_order1}
\end{equation} 
The part of this equation involving $N$ can be computed directly using the wave function $u(x,t)$ at time $t$.  To compute the exponential involving $D$ we use the fact that in the frequency domain, the partial derivative operator $\frac{\partial}{\partial x}$ is converted into $ik$, where $k$ is the frequency  associated with the Fourier transform. Then we take the Fourier transform of $u(x,t)$ recover the associate wave number, and compute 
$$
e^{-\frac{1}{2}idtk^{2}}\mathbb{F}\left[u(x,t)\right],
$$
where $\mathbb{F}$ denotes the Fourier transform. Then we take the inverse Fourier transform of the expression to find the solution in the physical space, yielding the final expression
$$
u(x, t+dt) = e^{idtN}\mathbb{F}^{-1}\left(e^{-\frac{1}{2} idtk^{2}} \mathbb{F} \left[u(x,t)\right]\right).
$$
We apply SSFM to Eq. \ref{NLSE2} where the nonlinear operator $N = |u|^{2}$ and the linear operator $L = \frac{1}{2}\frac{\partial^{2} }{\partial x^{2}} + \epsilon \delta(x)$.
In our numerical simulations we use the high-order SSFM, such as the Strang splitting based on: 
\begin{eqnarray*}
e^{i\left(L+N\right)\Delta t} & = & e^{i L\frac{\Delta t}{2}}e^{i N \Delta t}e^{i L\frac{\Delta t}{2}}+\mathcal{O}\big(\Delta t^3([L,[L,N]]+[N,[N,L]])\big).
\end{eqnarray*}
where $[L,N]=LN-NL$ denotes the commutator between $L$ and $N$. Thus, from $t$ to $t+\Delta t$
\begin{eqnarray}
u(x,t+\Delta t) & = & e^{i\left(L+N\right)\Delta t}u(x,t), \nonumber \\
                &\approx & e^{i L\frac{\Delta t}{2}}e^{i N \Delta t}e^{i L\frac{\Delta t}{2}}u(x,t).
                \label{split_order2}
\end{eqnarray}

\section{gPC collocation method}
We solve Eq. \ref{GP1} with the initial condition given by  Eq. \ref{initial1}  for both small and large values of $\epsilon$ by the gPC collocation method.
We use the gPC method for the solution of the NLSE using the Wiener-Askey scheme \cite{XiuBook, Xiu3}, in which Hermite, Legendre, Laguerre, Jacobi and generalized Laguerre orthogonal polynomials are used for modeling the effect of continuous random variables described by the normal, uniform, exponential, beta and gamma probability distribution functions (PDFs), respectively \cite{Chakraborty_Jung, Xiu2}.  These orthogonal polynomials are optimal for those PDFs since the weight function in the inner product and its support range correspond to the PDFs for those continuous distributions.

Following the standard gPC expansion, we assume that $u(x,t, \xi)$ is sufficiently smooth  in $\xi$ and has a converging expansion of the form 
$$u(x,t,\xi) = \sum_{k=0}^{\infty}\hat{u}_{k}(x,t)P_{k}(\xi),$$
where the orthonormal polynomials $P_{k}(\xi)$ correspond to the PDF of the random variable $\xi$ and satisfy the following orthogonality relation:
$$
\mathbf{E}[P_{k}P_{l}] := \int P_{k}(\xi)P_{l}(\xi)\rho(\xi)d\xi = \delta_{kl}.
$$
Here $\delta_{kl}$ is the Kronecker delta and $\rho(\xi)$ is the weight function. Note that the polynomials are normalized.

For the stochastic collocational approach we approximate $\hat{u}_{k}(x,t)$ as,

\begin{equation}
\hat{u}_{k}(x,t) = \sum_{j=0}^{Q}u\left(x,t,p^{j}\right)P_{k}\left(p^{j}\right)\alpha_{j}, \; k = 0, \cdots, Q,
\end{equation}
where $Q+1$ is the total number of the collocation nodes. 
Here $\left\lbrace p^{j},\alpha^{j} \right\rbrace$ is a set of nodes and weights, where $p^{j}$ and $\alpha^{j}$ denote the $j$-th node and its associated weights, respectively, in the random space $\Gamma$ such that
\begin{equation}
\mathbb{W}^{Q}\left[f\right]\equiv \sum_{j=0}^{Q}f\left(p^{j}\right)\alpha^{j},
\end{equation}
is an approximation of the integral
\begin{equation}
\mathit{I}\left[f\right]\equiv \int_{\Gamma}f(p)\rho(p)dp =  \mathbb{E}\left[f(p)\right],
\end{equation}
for sufficiently smooth functions $f(p)$, i.e,
$$
\mathbb{W}^{Q}\left[f\right] \rightarrow \mathit{I}\left[f\right], \; Q \rightarrow \infty.
$$

\begin{figure}
	\centering
		\includegraphics[width=0.8\textwidth]{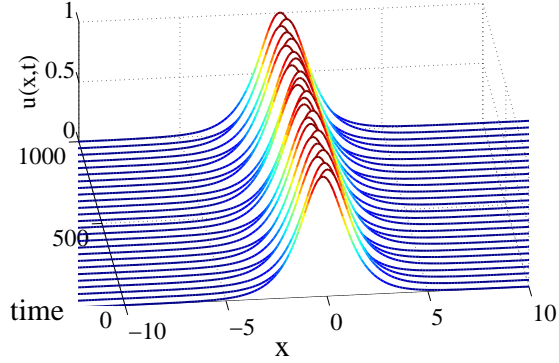}
		\caption{Soliton interaction with the defect with small strength $\left(\epsilon = 0.1\right)$, where initial velocity of the soliton is zero. The soliton is trapped and an oscillatory movement is observed.}
	 \label{fig:soliton0}
\end{figure}

In this paper we consider $V$ as the stochastic variable and we choose a collocation nodal set $\left\lbrace V^{j}, \alpha ^{j}\right\rbrace_{j=0}^{Q}$ in space $\Gamma$, where $V^{j}$ are the $j$th collocation points and $\alpha^{j}$ the corresponding weights. 
For each $j = 0,\cdots, Q$, we solve the problem given by Eqs. \ref{GP1} and \ref{initial1} with the parameters $\epsilon$ and $V^{j}$ and let the solution set be $\left\lbrace u_{0},\cdots , u_{Q}\right\rbrace$ where $u_{j}$ is the solution for $V = V_{j}$. For solving this deterministic equation, we employ the high-order SSFM. The approximate gPC expansion coefficients are
$$
\hat{u}_{m}(x,t) = \sum_{j=0}^{Q} u_{j}\left(x,t,V_{j}\right)\phi_{m}\left(V_{j}\right)\alpha_{j}, \; m = 0,\cdots , Q,
$$
where $\left\lbrace \phi_{m}\right\rbrace$ are the orthonormal polynomials. And finally we construct the $Q$th order gPC approximation
$$
u(x,t; V) \approx \sum_{m=0}^{Q}\hat{u}_{m}\left(x,t\right)\phi_{m}\left(V\right), \; \mathrm{where}\; V = \left\lbrace V_{0}, V_{2},  \cdots , V_{Q}\right\rbrace.
$$
\section{gPC collocation algorithm for solving NLSE}
\begin{figure}
	\centering
		\includegraphics[width=0.6\textwidth]{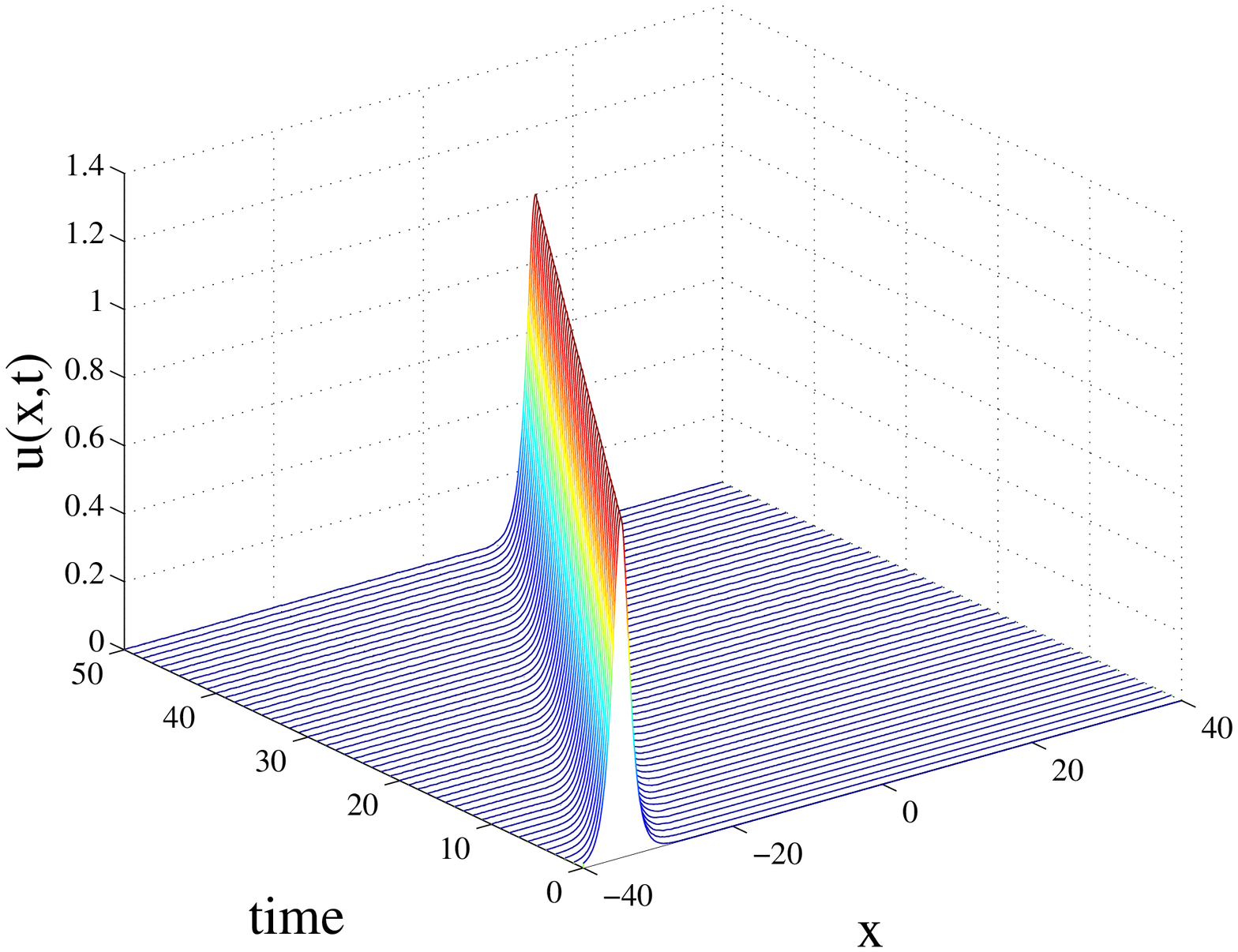}
		\includegraphics[width=0.6\textwidth]{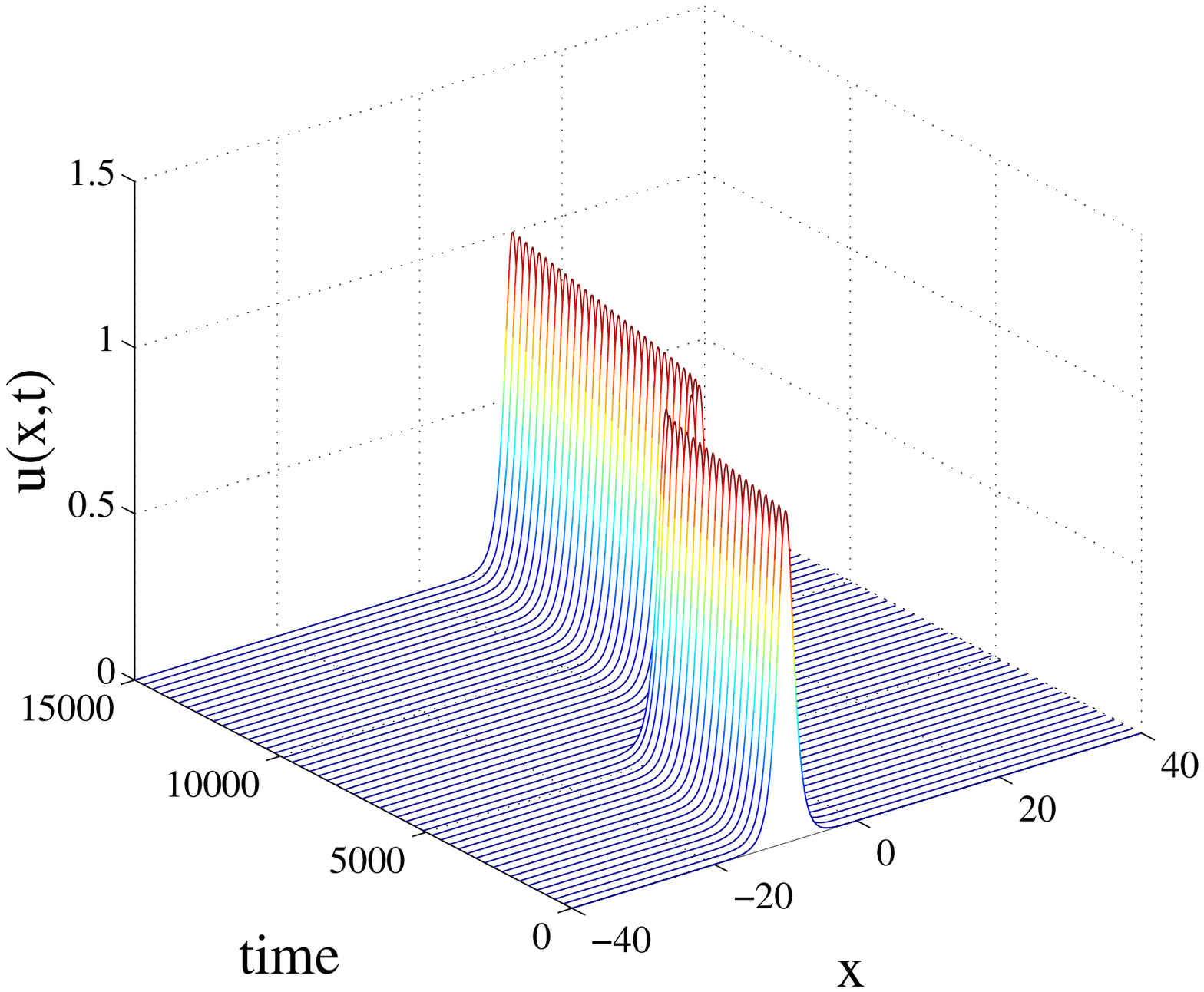}
		\includegraphics[width=0.60\textwidth]{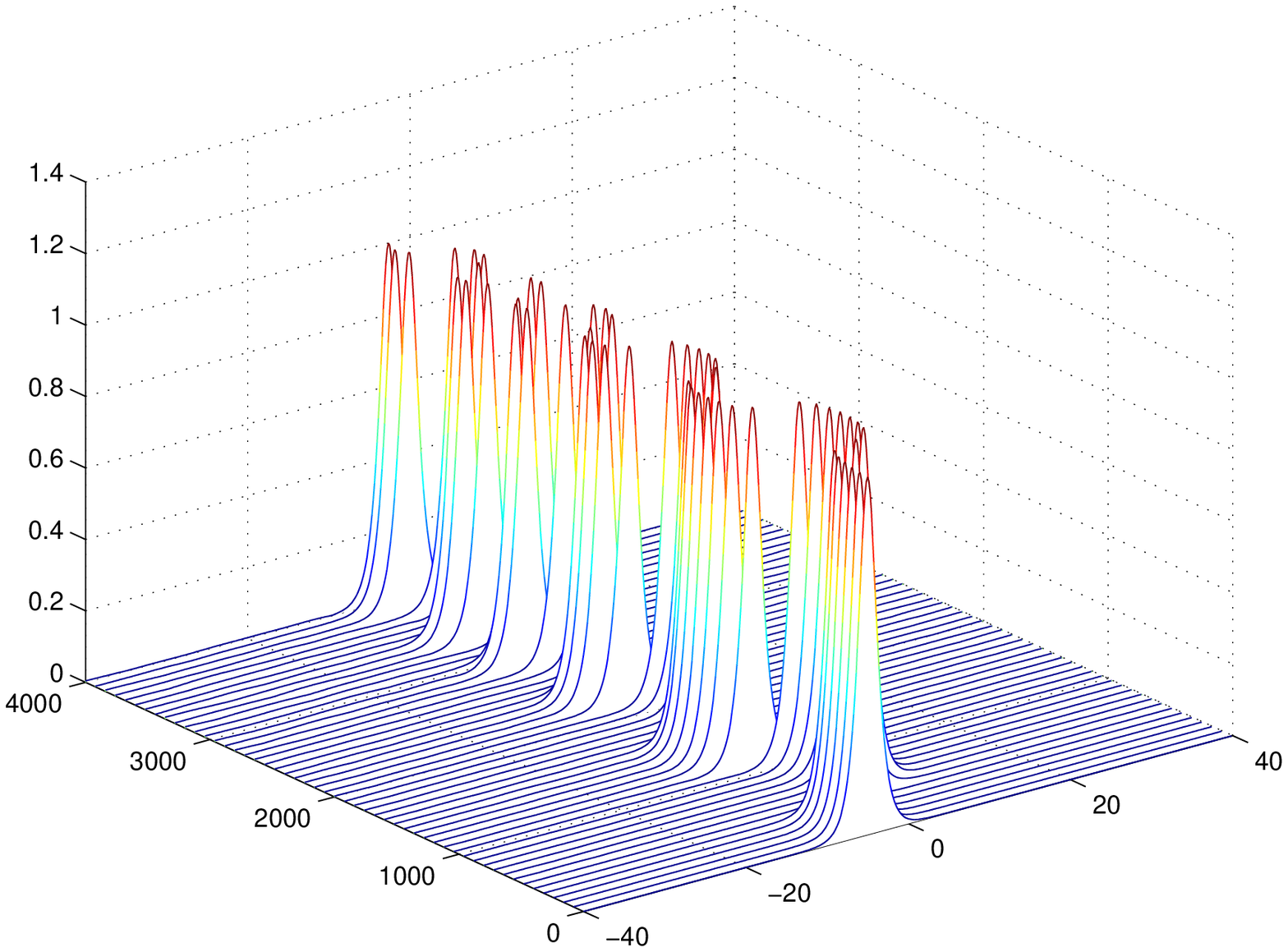}
		\caption{Top: Soliton solution without any defect. Middle: Soliton passes through the defect with the initial velocity $V = 0.001$ and the defect amplitude $\epsilon = 0.1$. Bottom: Soliton trapped by the defect with the initial velocity $V = 0.003$ and the defect amplitude $\epsilon = 0.5$}
	 \label{fig:soliton1}
\end{figure}
The following algorithm describes how to calculate the critical velocity by using the gPC collocation method. 

We use the gPC method  to find the critical velocity $V_{c}$ efficiently for any given $\epsilon$. Here the soliton velocity $V$ is the stochastic variable. Suppose we know in advance that the critical velocity $V_{c}$ lies between $V_{a}$ and $V_{b} \; \left(V_{a} < V_{b}\right)$ and consider $V$ has a uniform distribution over $\left[V_{a}, \; V_{b}\right]$. Since the distribution is uniform, we use the Legendre polynomials for expanding the solution in the random space.  For this purpose we choose $N+1$ Gauss-Legendre quadrature points with the weights. Let the set $\left\lbrace \alpha_{i},\;\omega_{i}\right\rbrace_{i=0}^{N}$ describe the $\left(N+1\right)$ quadrature points $\alpha_{i}$ and the corresponding weights $\omega_{i}$.

Now find the solution of Eq. \ref{GP1} for each $V = \alpha_{i}$ by using the high-order SSFM. For this purpose one must use a sufficiently large computational domain and sufficiently long time interval. We set up the domain size and the computational time in such a way that no solution leaves  the domain yet with the given final time. For example when $\epsilon = 0.3$, we use the domain size $\left[-L, \;L\right]= \left[-40, \; 40 \right]$ and the final time $t_{f} = 12000$. We are solving the NLSE for $u_{j}\left(x, t, V_{j}\right)$ for all $V_{j}$ with the same final time.

We reconstruct the soliton solution for each simulation for  $ x \in \left[-L, \; L^{'}\right]$ at the final time. $L^{'}$ is chosen in such a way that only the trapped solutions exist inside $\left[-L, \; L^{'}\right]$. We know if the solution is trapped, it would stay around the position of the defect (in our case at $x = 0$). So $L^{'}$ must be close to zero. In our computation, we choose $L^{'}$ where the mean solution vanishes near $x = 0^{+}$. For the $i$-th quadrature point $\alpha_{i}$, we denote the solution by $u_{i}\left(x,T_{f}, \alpha_{i}\right)$.

\begin{figure}
	\centering
		\includegraphics[width=0.47\textwidth]{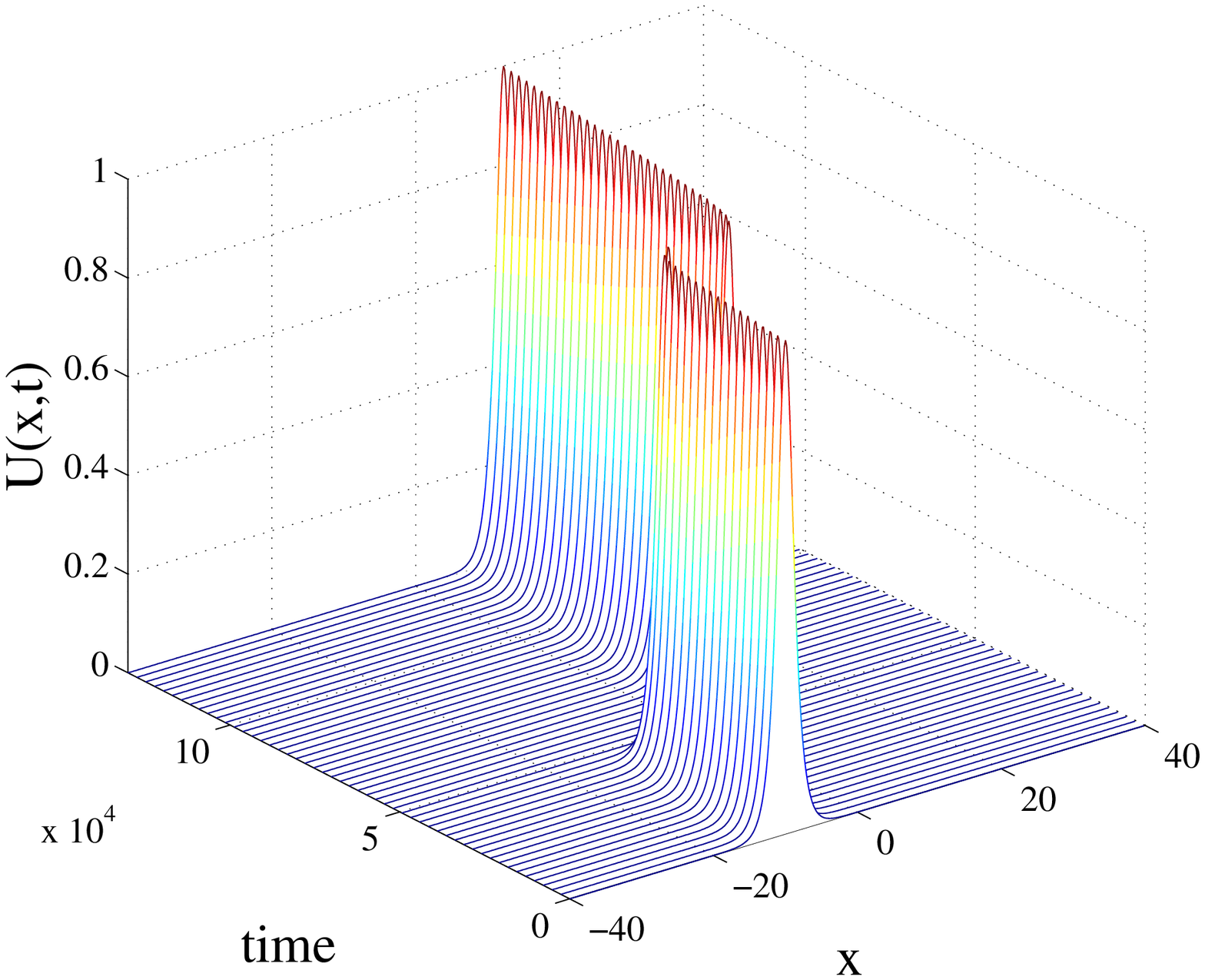}
		\includegraphics[width=0.47\textwidth]{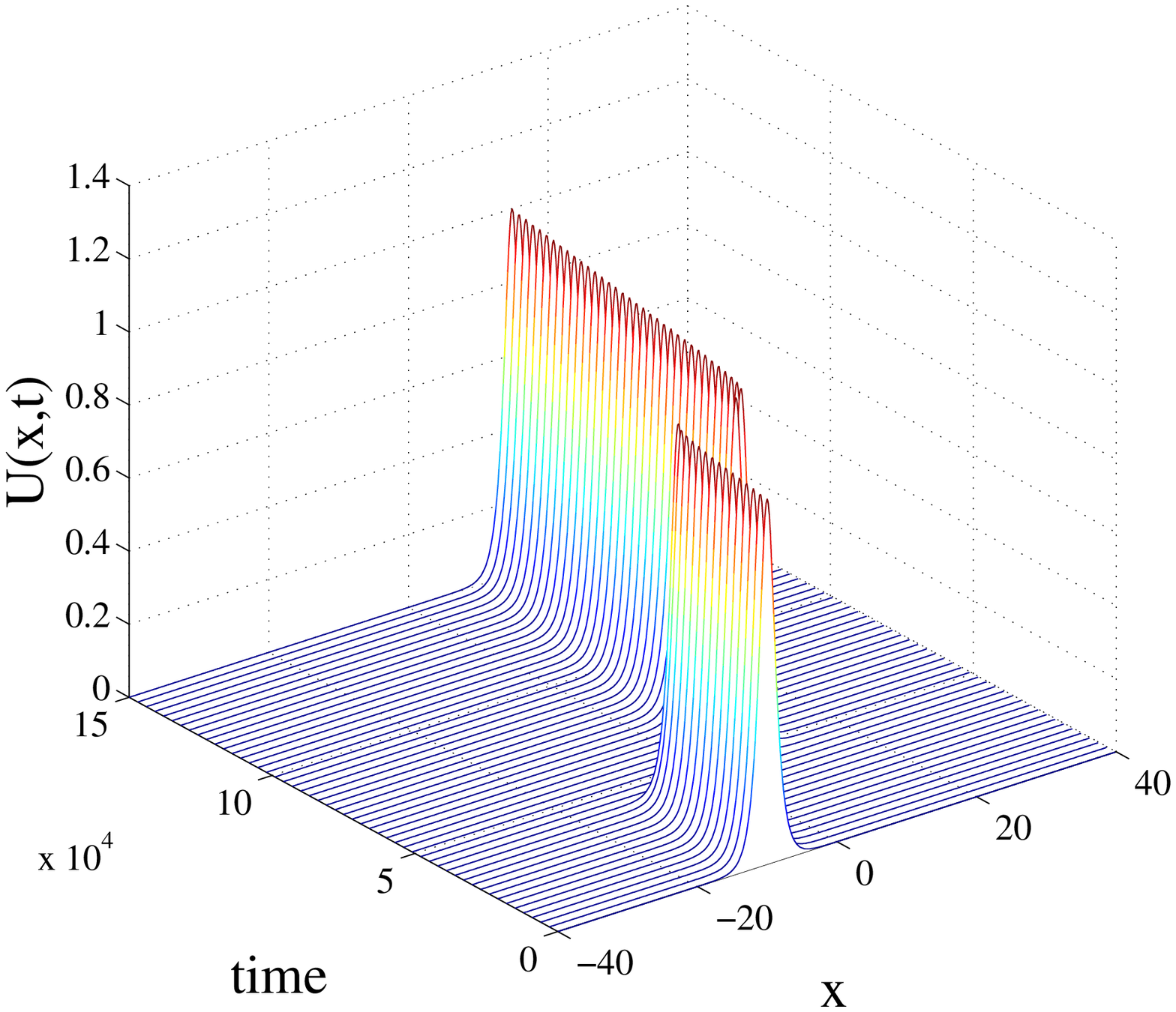}
		\includegraphics[width=0.47\textwidth]{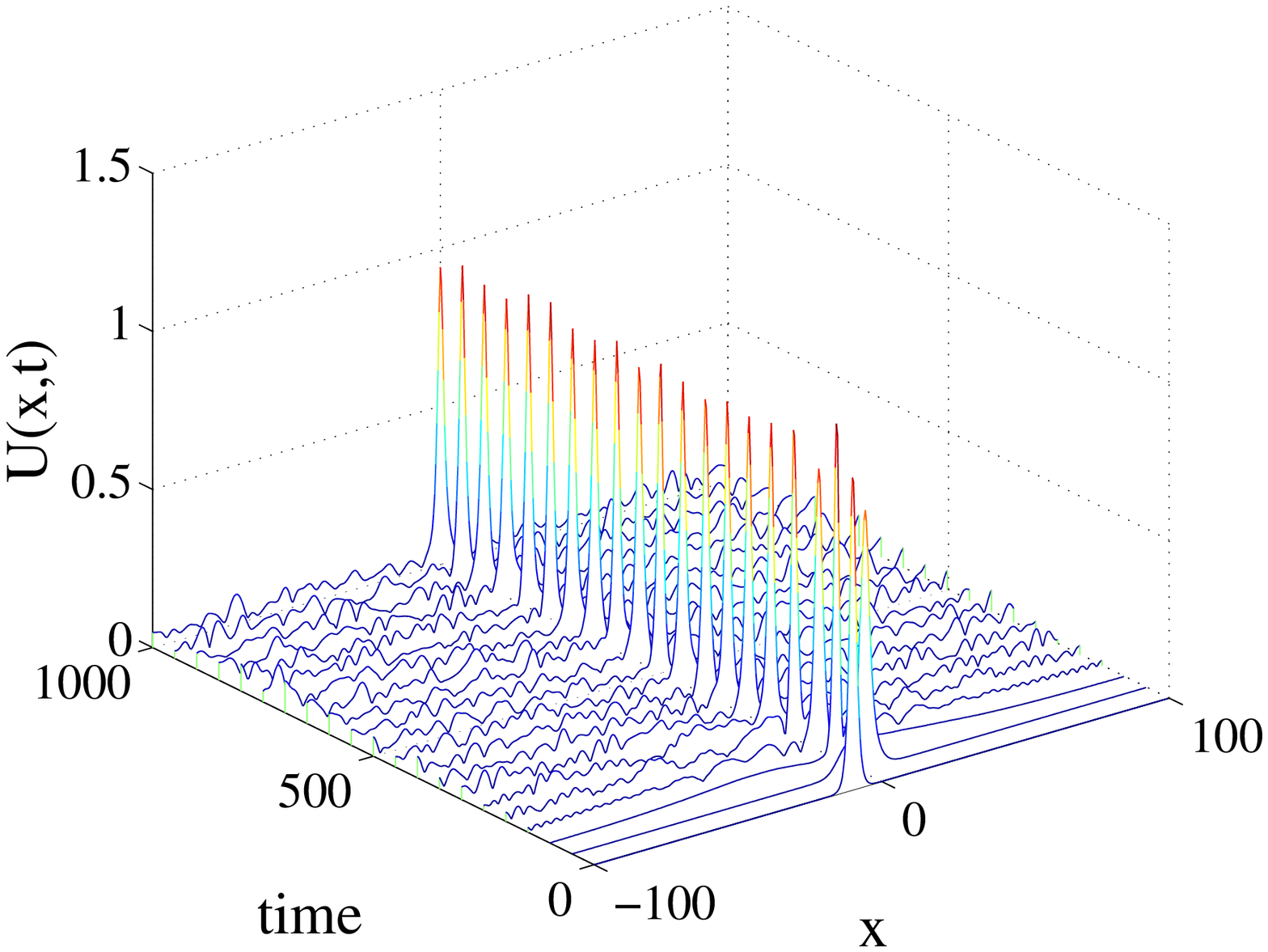}
		\includegraphics[width=0.47\textwidth]{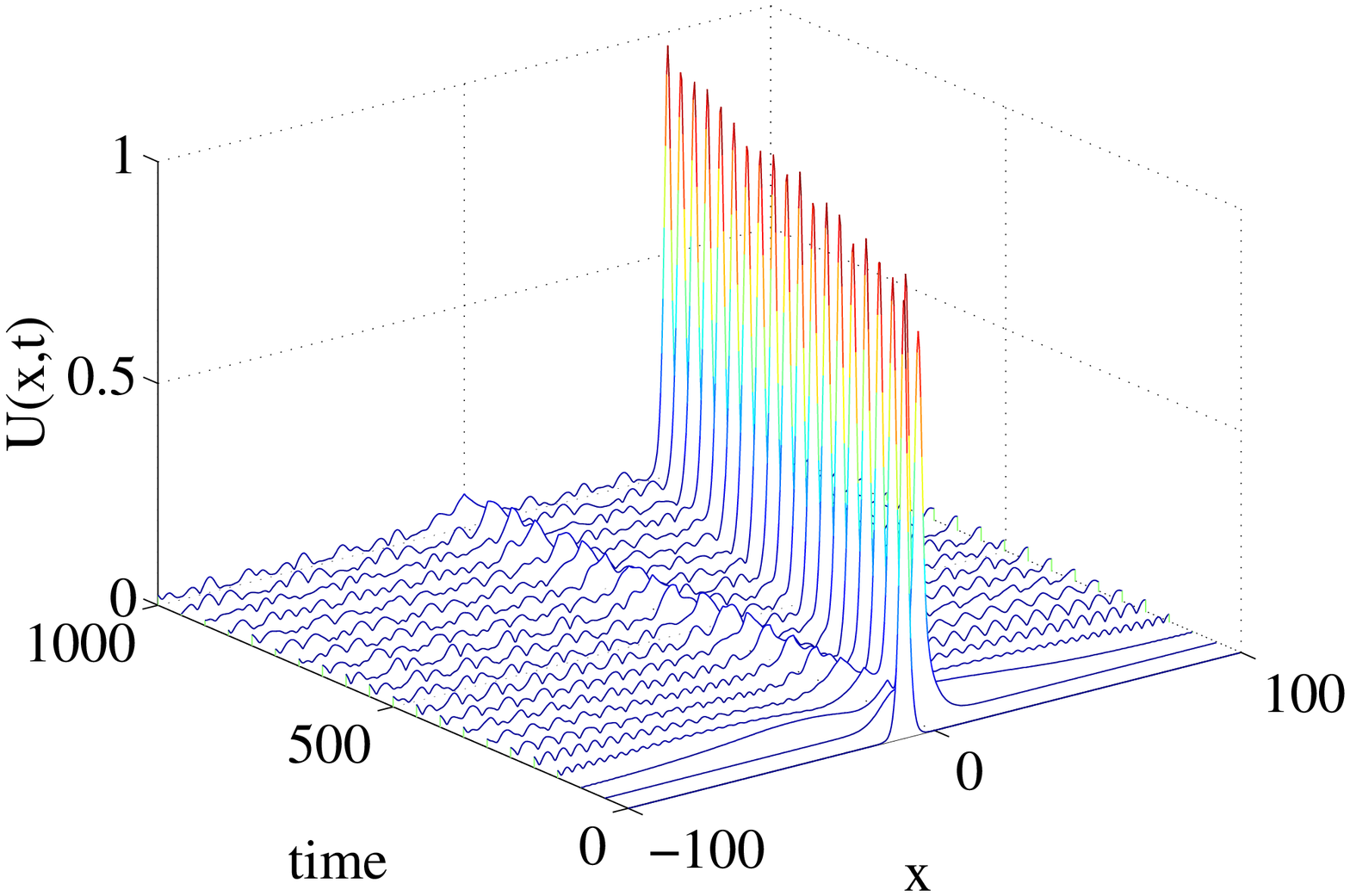}
	    
		\caption{Top left: Soliton transmitted through the defect when $\epsilon = 0.08$ and $V = 5 \times 10^{-5}$. Top right: Soliton transmitted through the defect when $\epsilon = 0.1$ and $V = 8 \times 10^{-5}$. Bottom left: Soliton is trapped by the defect when $\epsilon = 4.5$ and $V = 0.220048$.
		Bottom right: Soliton is transmitted through the defect when $\epsilon = 4.5$ and $V = 0.23995187$. For $\epsilon = 4.5$, radiation effect is clearly visible. }
	 \label{fig:soliton3}
\end{figure}
\begin{figure}
	\centering
		\includegraphics[width=0.60\textwidth]{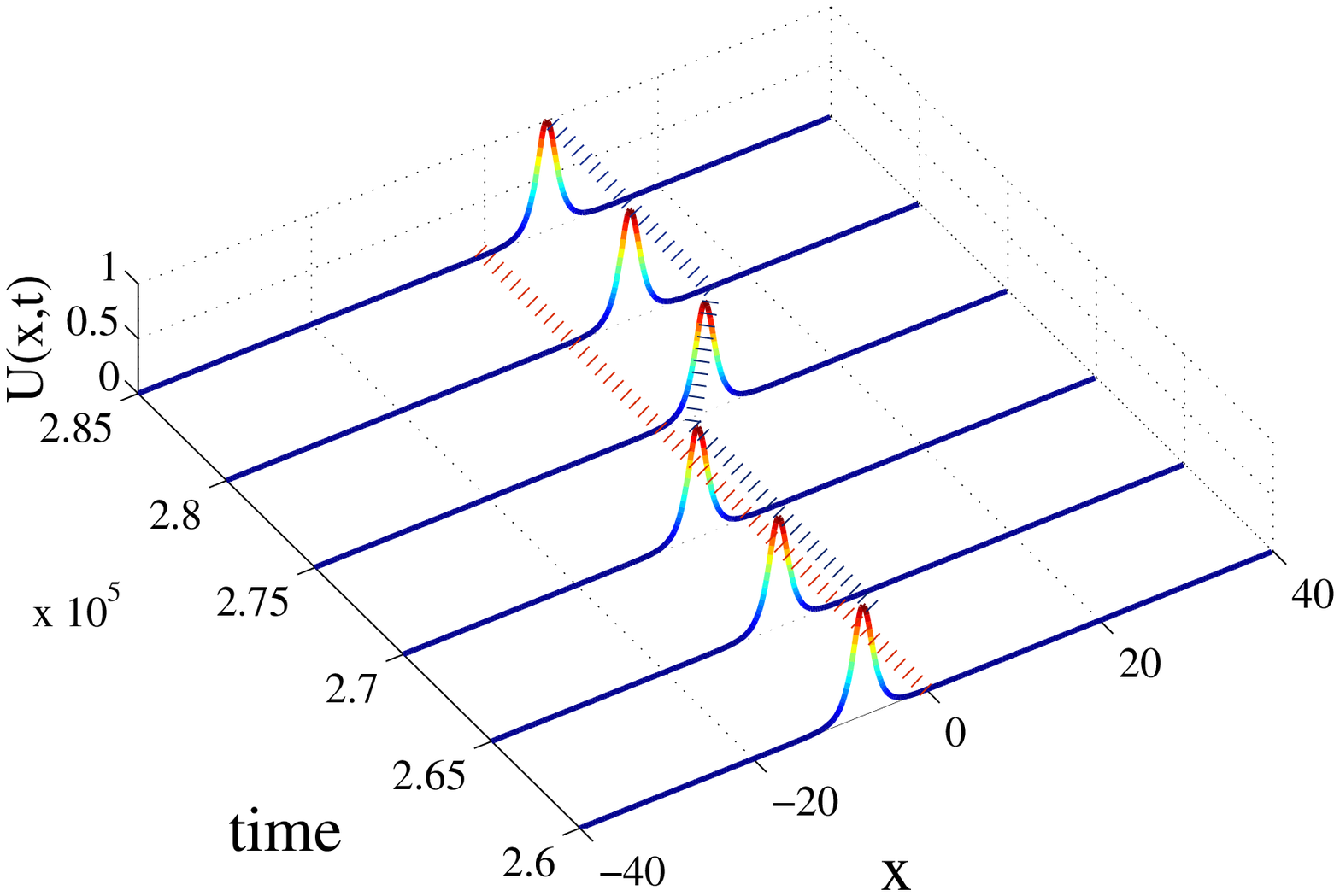}
		\includegraphics[width=0.60\textwidth]{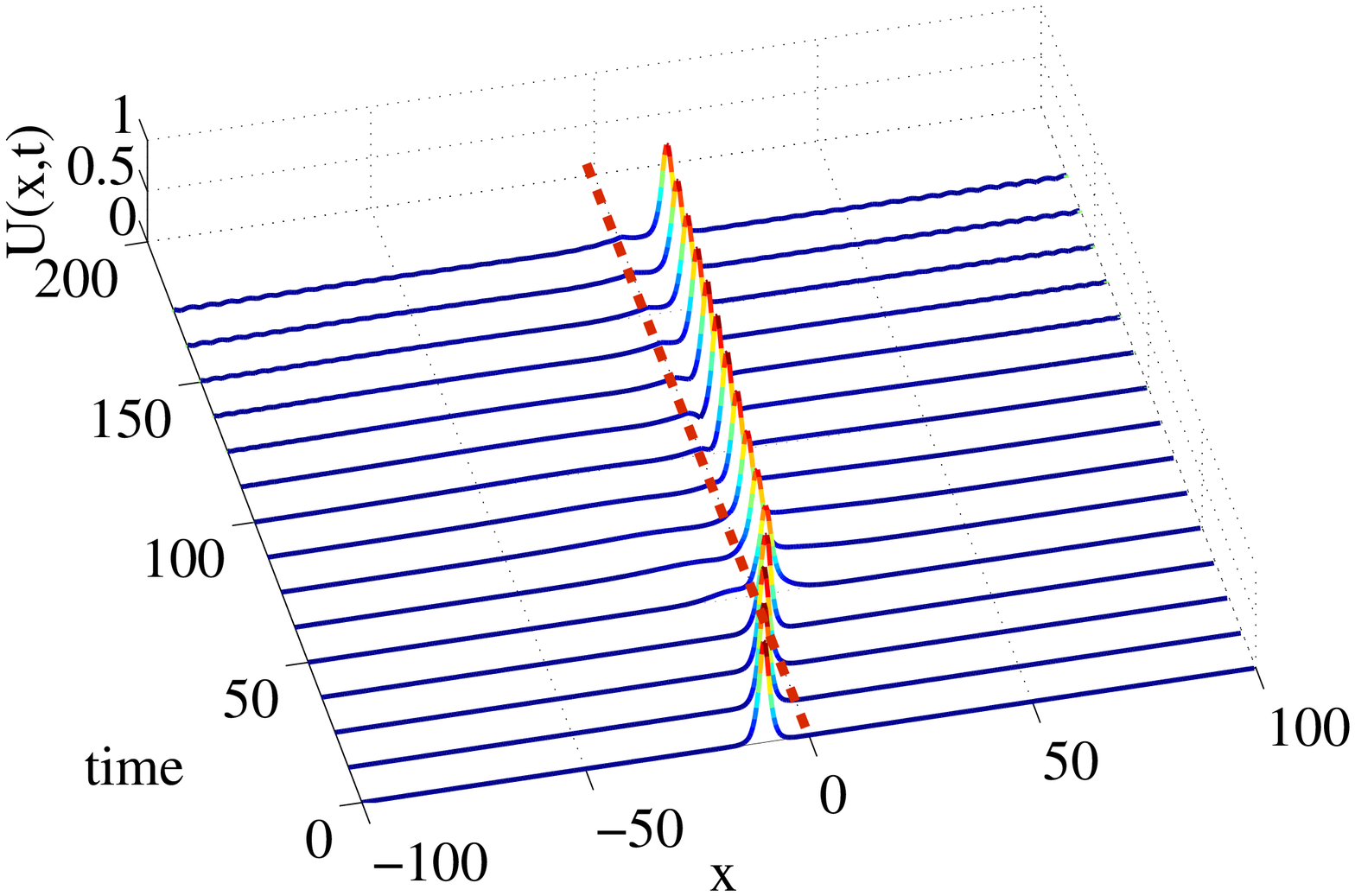}
	    \includegraphics[width=0.60\textwidth]{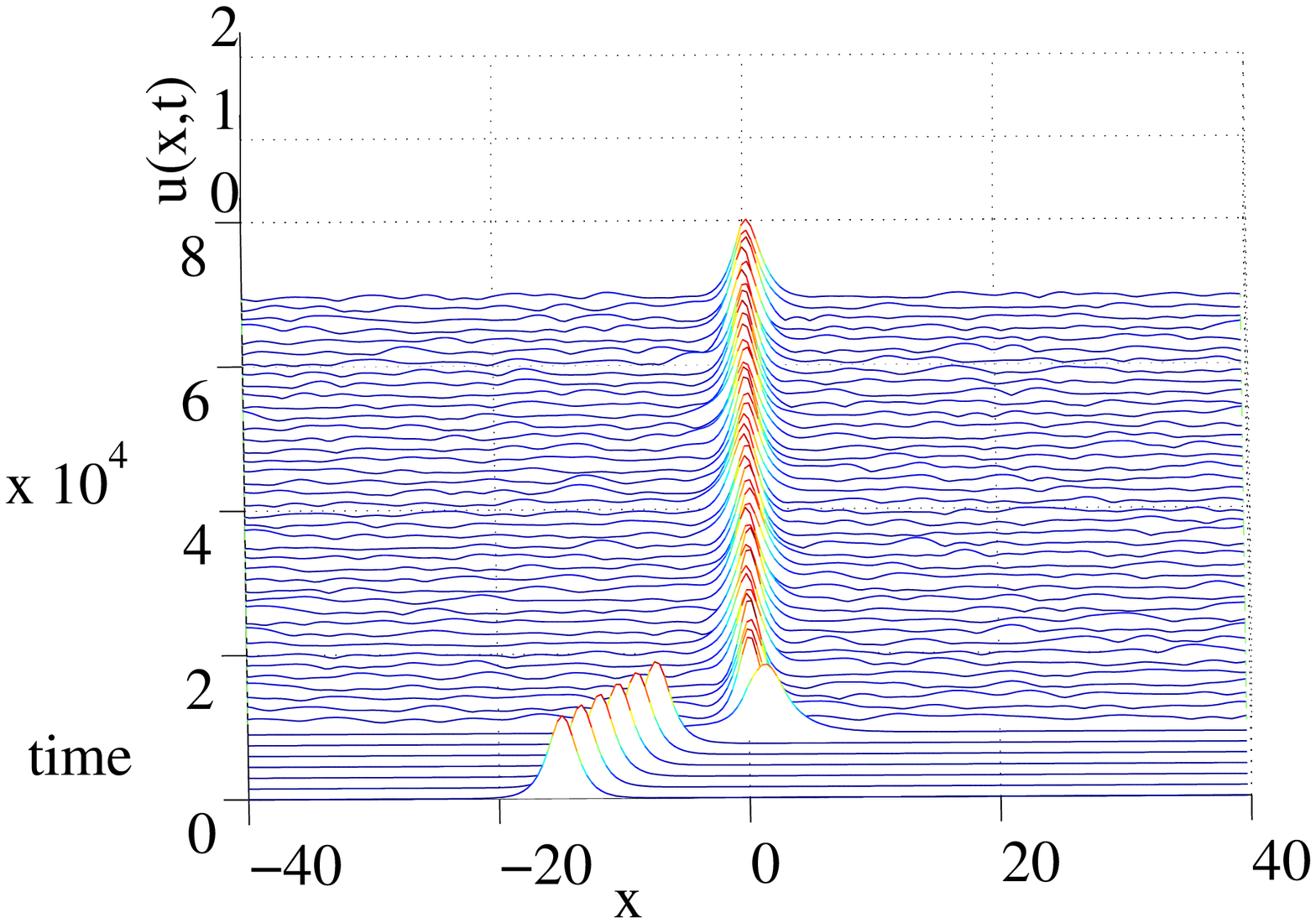}
		\caption{The nonlinear interaction of the soliton with the defect (the dotted line). Top: The nonlinear interaction is prominent when the soliton  hits the deffect $\left(\epsilon = 0.3\right)$ with small velocity $\left(\sim 10^{-5}\right)$   compared to the  interaction with the high velocity $\left(\sim 10^{-1}\right)$ where  $\epsilon = 4.5$ (middle). Bottom: The interaction of slowly moving soliton $\left(V \sim 10^{-5}\right)$) with the defect with high value of $\epsilon \left( = 4.5\right)$.}
	 \label{fig:poincare}
\end{figure}
Evaluate the approximate gPC expansion coefficients by
$$
\hat{u}_{m}\left(x, T_{f}\right) = \sum_{i=1}^{Q}u_{i}\left(x,T_{f},\alpha_{i}\right) L_{m}\left(\alpha_{i}\right)\omega_{i},
$$
where $\left\lbrace L_{m}\right\rbrace_{m = 0}^{Q}$ is the set of  Legendre polynomials and $\omega_i$ are the quadrature weights. The full gPC solution is given by 
\begin{equation}
u(x,T_{f},\alpha) = \sum_{k=0}^{Q}\hat{u}_{k}\left(x,T_{f}\right)L_{k}(\alpha).
\label{sol_gpc1}
\end{equation}
The mean solution is given by the 1st mode \cite{Xiu2}, i.e
\begin{equation}
\hat{u}_{0}\left(x, T_{f}\right) = \sum_{i=0}^{Q}u_{i}\left(x,T_{f},\alpha_{i}\right)L_{0}\left(\alpha_{i}\right)\omega_{i}.
\label{sol_gpc_mean}
\end{equation} 
From Eq. \ref{sol_gpc_mean}, one can construct the average energy $\bar{E}$ of the system between $\left[-L, \; L^{'}\right]$ and $\left[-L, \; L\right]$ at the final time, that is, 
\begin{eqnarray}
\bar{E}_{L} =\frac{1}{2} \int_{-L}^{L}\left|\hat{u}_{0}\left(x, T_{f}\right)\right|^{2}dx,\; \bar{E}_{L^{'}} =\frac{1}{2} \int_{-L}^{L^{'}}\left|\hat{u}_{0}\left(x, T_{f}\right)\right|^{2}dx.
\label{definition}
\end{eqnarray}
Suppose that among $N$ solutions, $N_{1}$ solutions are trapped inside $\left[-L, \;L^{'} \right]$. 
Then $N_{1}$ can be estimated for large $N \rightarrow \infty$ by 
\begin{eqnarray}
\frac{N_{1}}{N} = \frac{\bar{E}_{L^{'}}}{\bar{E}_{L}}, \quad  \frac{N_{1}}{N} = \frac{V_{c}-V_{a}}{V_{b}-V_{a}}, \nonumber 
\end{eqnarray}
where $V_{c}$ is the critical velocity for given $\epsilon.$ So $V_{c}$ is evaluated by
\begin{equation}
V_{c}= V_{a} + \left(V_{b} - V_{a}\right) \frac{\bar{E}_{L^{'}}}{\bar{E}_{L}}.
\label{v_critical}
\end{equation}
If we increase the number of quadrature points, then the critical velocity can be determined more accurately. For our simulations we used $24$  Gauss Legendre quadrature points and obtained spectral accuracy of $\sim 10^{-12}$.
Figure \ref{fig:soliton6} shows  the spectral convergence of the error of the critical velocities  with the increasing number of the quadrature points. 


\noindent
\textbf{Remark:}

\noindent
\textit{The solution $u(x,t,V)$ has possibly a jump at $V = V_c$  for $t \rightarrow \infty$ because of the critical behavior of the soliton solution around the potential. This means that the spectral reconstruction of $u(x,t,V)$ for any $V \in [V_a, V_b]$ using ${\hat u}_l (x,t), l = 0, \cdots, Q$ may fail to converge to the right solution due to the discontinuity at $V = V_c$. This was also addressed in our previous work for the critical behavior of the soliton solution for the sine-Gordon equation \cite{Chakraborty_Jung}.  Here note that the proposed method in this paper uses only the first moment ${\hat u}_0(x,t)$ to estimate the critical velocity $V_c$ but not the reconstruction of $u(x,t,V)$. The mean solution, ${\hat u}_0(x,t)$ is convergent. }

In Eq. \ref{v_critical}, the convergence of $V_c$ mainly depends on $R := \frac{\bar{E}_{L^{'}}}{\bar{E}_{L}}$. As the definition in Eq. \ref{definition}, the convergence of $R$ then depends on how ${\hat u}_0(x,t)$ converges with $N$.  In our previous work \cite{JungSong}, it was proven that ${\hat u}_0(x)$ converges fast enough although the original function $u(x,V)$ is discontinuous in the random variable $V$.  As we will discuss in the next section, numerical results in Section 5 (Figure 7) implies that $R$ shows spectral convergence with $N$. 

\begin{figure}
	\centering
		\includegraphics[width=0.47\textwidth]{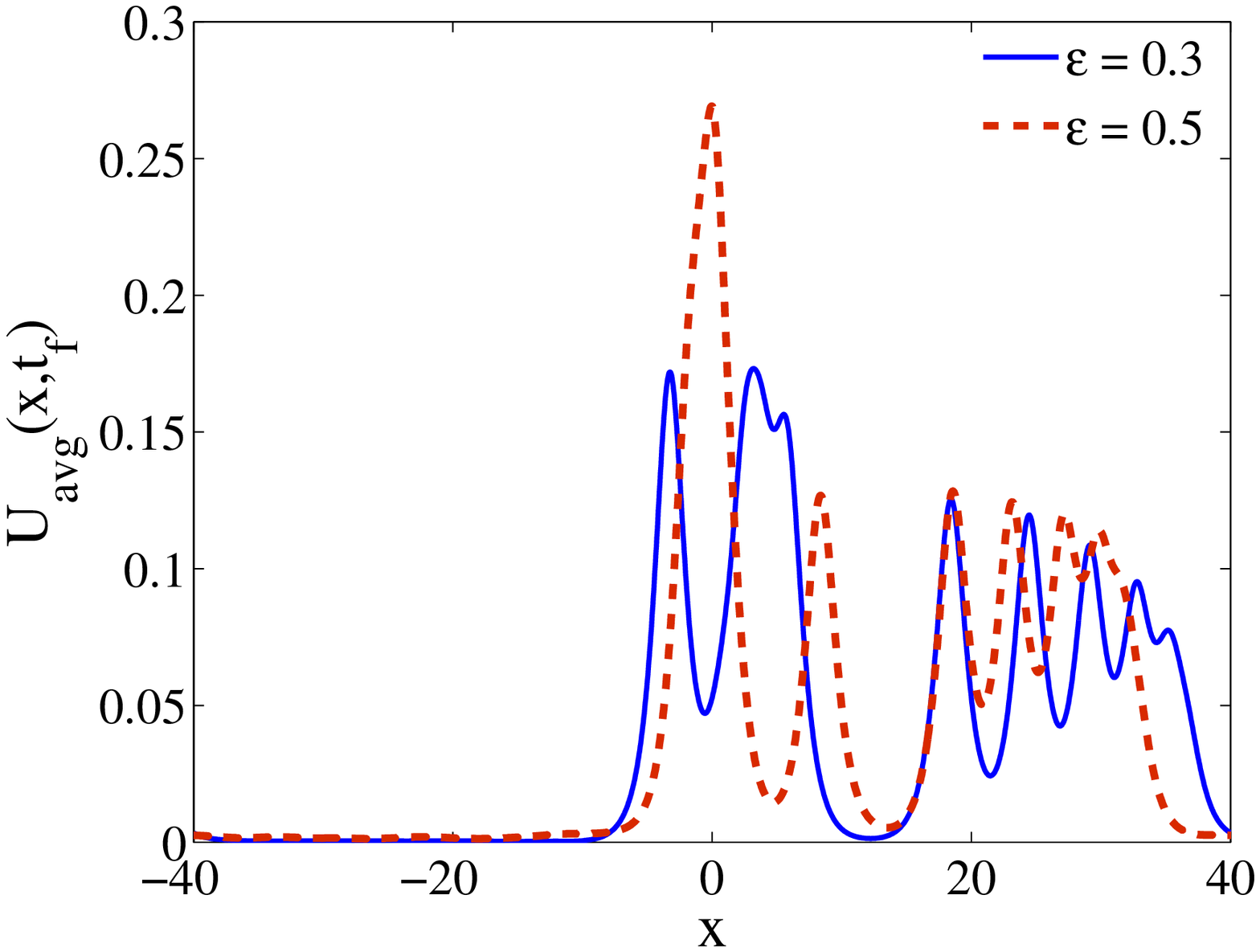}
		\includegraphics[width=0.47\textwidth]{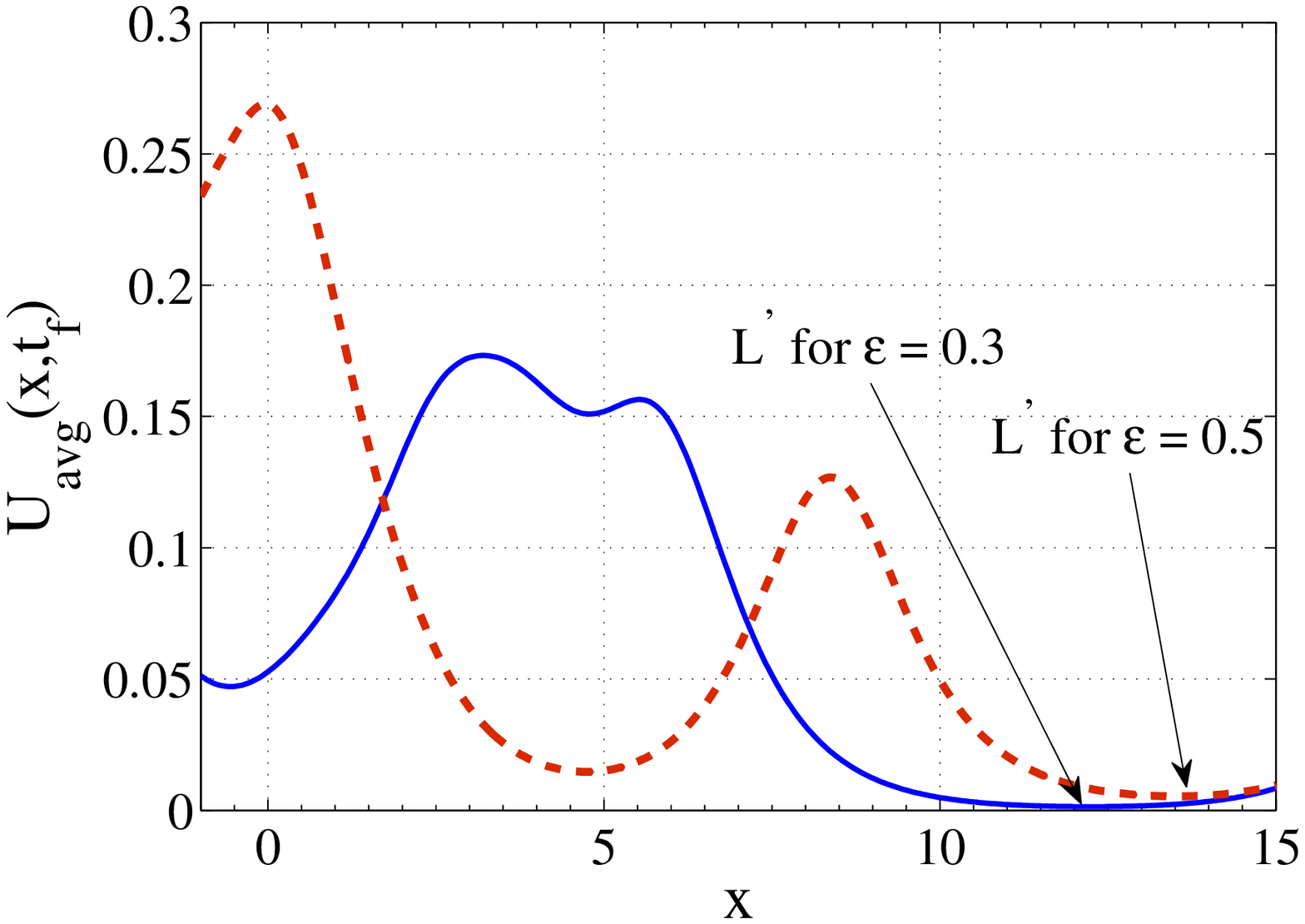}
		\includegraphics[width=0.47\textwidth]{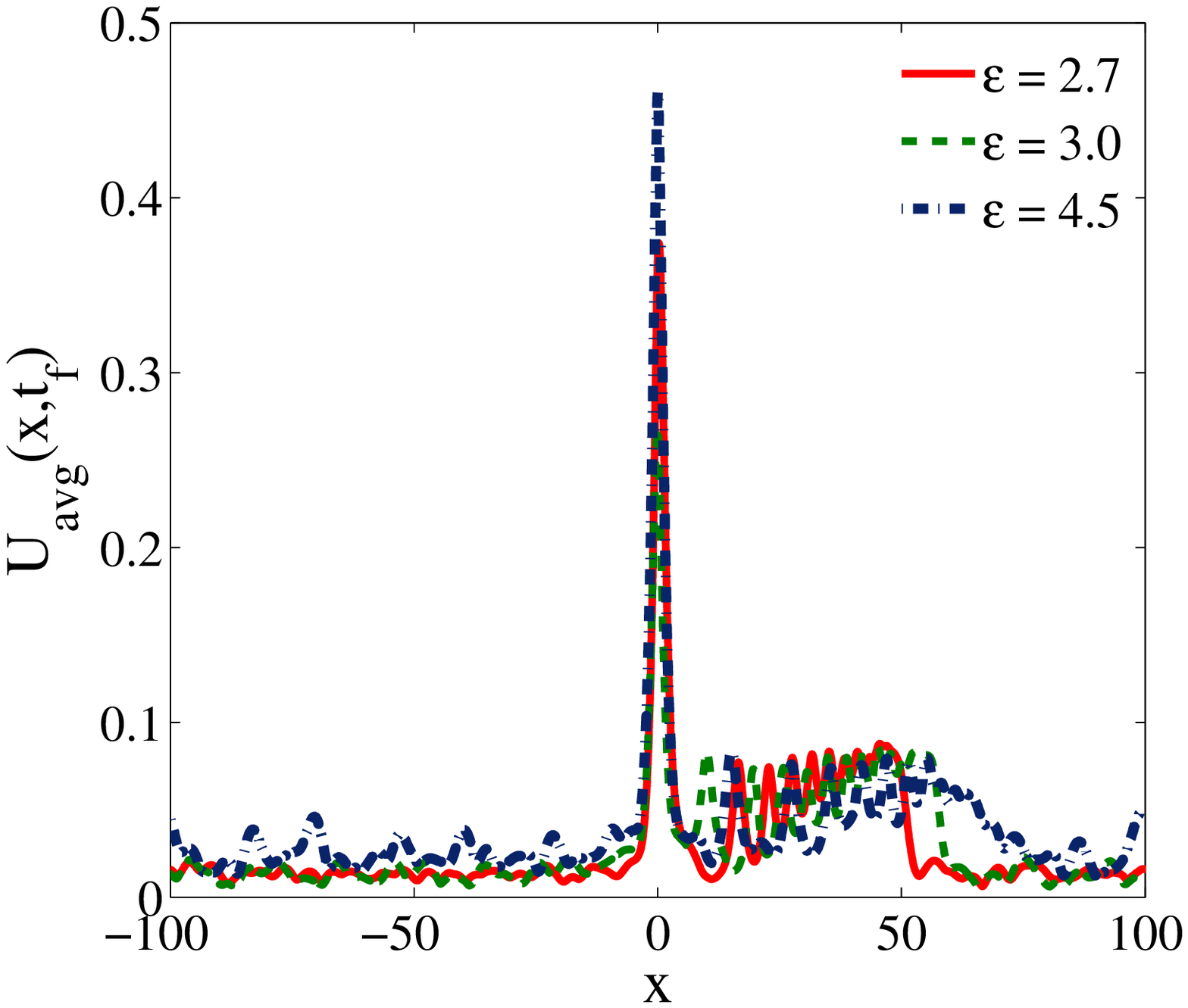}
		\includegraphics[width=0.47\textwidth]{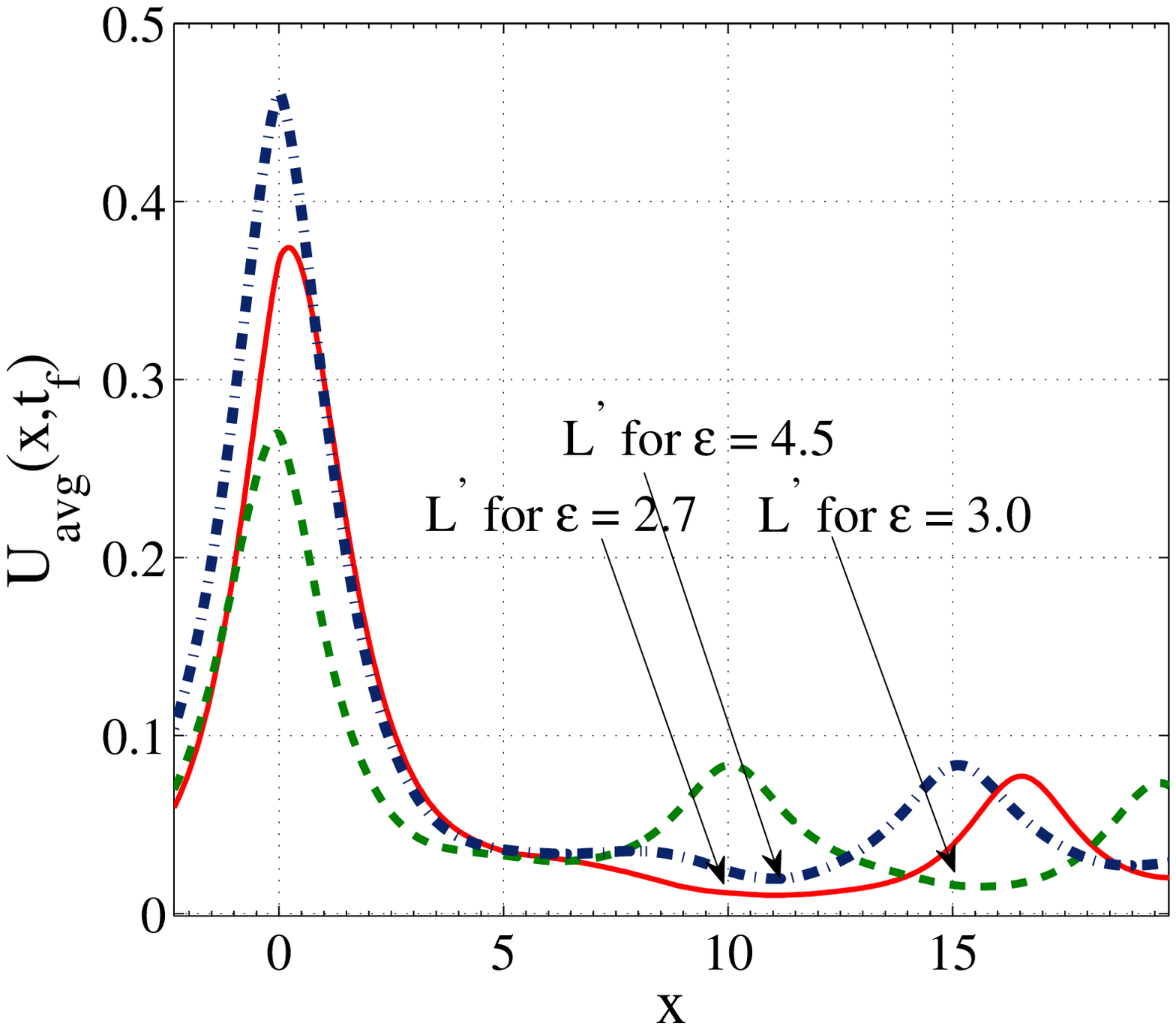}
		
		\includegraphics[width=0.55\textwidth]{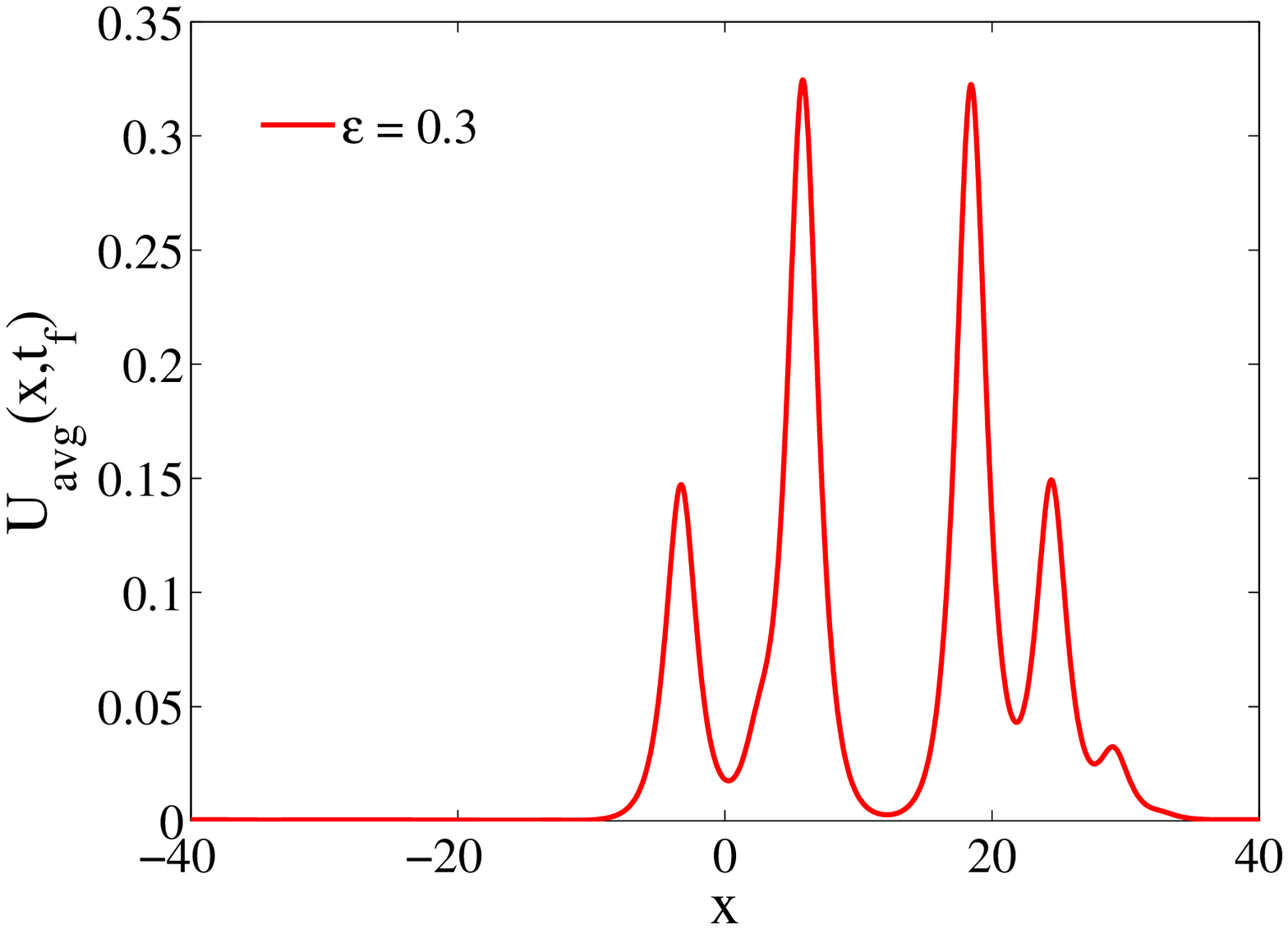}
	    
		\caption{ First mode (mean) of the gPC expansion for different $\epsilon$. Top and middle: The Legendre chaos. Bottom: The Hermite chaos. Right figures of the top and middle panels show the locations of $L^{'}$ for different $\epsilon$.  }
	 \label{fig:soliton4}
\end{figure}
\section{Numerical results}
\begin{figure}
	\centering
		\includegraphics[width=0.8\textwidth]{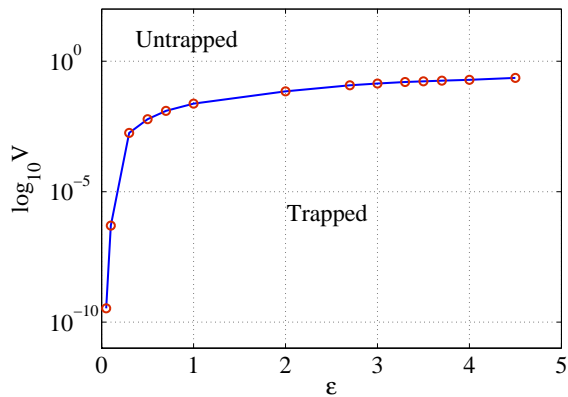}
	    
		\caption{ Critical Velocity vs. $\epsilon$ where $\epsilon \in \left[0.05, \; 4.5\right]$.  }
	 \label{fig:soliton5}
\end{figure}

We first consider the high value of $\epsilon$, say $\epsilon = 2.7$. By doing few Monte-Carlo simulations we roughly estimate the interval  $V \in \left[V_{a},\; V_{b}\right], \; V_{c} \in \left[V_{a},\; V_{b}\right] $ where $V_{c}$, the critical velocity may be located. For $\epsilon = 2.7$, we use $V_{a} = 0.1$ and $V_{b} = 0.14$. Since for moderate and high values of $\epsilon$, the simulation time is relatively less than the simulation time with smaller range of $\epsilon$, we follow the same procedure to find the suitable intervals. But for the small value of $\epsilon$, i.e. $\epsilon < 1.0$, where the simulation time is long, we use the extrapolation of $V_{c}$ from the previous $\epsilon$ to get the rough estimate of the interval.

To apply the gPC collocation method, one also needs to find the value of $L^{'}$. We do not have any fixed $L^{'}$ which can serve for all $\epsilon$. Instead, we have different $L^{'}$ for different $\epsilon$. A heuristic approach is used to find $L^{'}$. For  the given value of $\epsilon$, we construct the mean solution by Eq. \ref{sol_gpc_mean}. Since some solutions are trapped and some of them are transmitted, there are few bumps near the defect and few bumps are far  from the defect. Clearly there exists a separation point between these two groups of bumps. Ideally the $x-$coordinate of this point would be zero but due to the domain truncation, radiation effect etc. it may not be equal to zero.  By observing the graph carefully we can easily find the separation point which we use as $L^{'}$. For the small and moderate values of $\epsilon$, determining accurately $L^{'}$ is easy, but for the high values of $\epsilon$, we need extra care. For the high value of $\epsilon$, the values of $V_{a}$, $V_{b}$ are also high and we can not run the simulations for a long time because some solutions may leave the domain and re-enter the domain from the other side due to the periodic boundary conditions. So in this case we need to study the bumps carefully to locate $L^{'}$. In Figure \ref{fig:soliton4}, the zoomed graphs of the mean solution of each $\epsilon$ are given in the right panel of the top and middle figures. We find that for $\epsilon = 0.3, \; L^{'} = 12$ and for $\epsilon = 0.5, \; L^{'} = 13.5$. Similarly for $\epsilon = 2.7, \; L^{'} = 10$ and for $\epsilon = 3.0, \; L^{'} = 15$.

Figure \ref{fig:soliton0} presents the interaction of the soliton with the $\delta$-function. Here we choose the initial velocity, $V_{0} = 0$ and the potential strength $\epsilon = 0.1$. The soliton is located at $x_{0} =  -0.3$ initially, which is inside the influence zone of the potential. The nonlinear interaction is observed and the soliton solution exhibits an oscillatory behavior along the line $x = 0$. This case was discussed in \cite{Holmer1, Holmer2}. But such an  initial condition may not necessarily satisfy the given equation. The initial position of the soliton must be out of the influence zone of the potential and the soliton must be allowed to move freely before it hits the defect.
In all cases we consider the starting point of the soliton $\left(x_{0}\right)$ is far from the position of $\delta$-function, i.e. outside the influence region of the potential.
Figure \ref{fig:soliton1} shows the behavior of the soliton solutions in three different cases.  When $\epsilon = 0$, that is the case when there is no $\delta$-function, the soliton solution passes unperturbedly. But for nonzero $\epsilon$, the soliton behaviour depends on its initial velocity. For $\epsilon = 0.1$, the soliton passes through the defect for $V = 0.001$ and for $\epsilon = 0.5$ and $V = 0.003$, soliton is trapped by the defect. For both cases, the soliton passed or trapped as a whole. There is no radiation due to the small soliton velocities \cite{Cao_Malomed}.

Figure \ref{fig:soliton3} represents the long time simulations for $\left(\epsilon,\;V\right) = \left(0.08,\; 5\times 10^{-5}\right)$ (top left), $\left(0.1,\; 8\times 10^{-5}\right)$ (top right), $\left(4.5,\; 0.220048\right)$ (bottom left) and $\left(4.5,\; 0.23995187\right)$ (bottom right).
For the case that $\epsilon$ is small and $V$ is also very small accordingly, the soliton is transmitted through the defect without any radiation. But for the high value of $\epsilon$, usually greater than $2.7$, where the critical velocity is also high, the radiation effect is observed due to the soliton-defect interaction. The bottom panel of  Figure \ref{fig:soliton3} exhibits the radiation effect for $\epsilon = 4.5$. For both the ``trapped" and ``transmitted" situations, the radiation effect is observed.  The bound state effect is also observed in the bottom right, the details of which was discussed in \cite{Holmer}.

Figure \ref{fig:poincare} shows the nonlinear interactions of the soliton with different soliton velocities. When the soliton velocity is small, nonlinear property dominates as shown in Figure \ref{fig:soliton3}. During the time of interaction with the defect (the dotted line), the soliton velocity increases and after crossing the defect, the velocity turns into its previous value. When the soliton velocity is high, the linear effect dominates and  the soliton velocity does not changes during the collision but the direction of the propagation changes. That is the soliton continues its motion with the same velocity.
When a slowly moving soliton hits the defect with high strength $\left(\epsilon = 4.5\right)$, the soliton is trapped by the defect but  due to  the nonlinear interactions, radiations and transmissions are also seen (bottom figure).

Figure \ref{fig:soliton4} shows the mean solutions at the final time. This is the first mode of the solution by the gPC collocation method. Here we used $V$ as a stochastic variable,  $V \in \left[V_{a}, \; V_{b} \right]$ and  $V_{a}$ and $V_{b}$ are different for different values of $\epsilon$. We used both the Legendre and Hermite chaos. We need to consider the uniform distribution and normal distribution for the Legendre and Hermite chaos respectively. In Figure \ref{fig:soliton4}, the  figures in the top panel are obtained using the Legendre chaos for $\epsilon = 0.3,\; 0.5$. 
Those solitons  that are trapped by the defect are confined around the position of the defect. In our case, the defect, the $\delta$-function is located at $x = 0$. There are multiple peaks in the mean solution, but around $x = 0$ the peaks are higher than the others, which implies that some solitons are trapped, and the rest are transmitted. These figures are used to locate the position of $L^{'}$. If we see the zoomed figure in the right panel, we easily locate $L^{'}$ for different $\epsilon$.

For the middle panel figures in Figure \ref{fig:soliton4}, we plotted the mean solutions and zoomed one for $\epsilon = 2.7,\; 3.0 \; \mathrm{and}\; 4.5$. The sharp peaks at $x = 0$ imply that the most of the solutions are trapped in that range of $V$ and  some of them are transmitted. We already mentioned that in this region of such a large value of $\epsilon$, the radiation effects are visible, which are also showed in the figure. The values of $L^{'}$ are pointed for different $\epsilon$ values in the figure.
Same explanation for $\epsilon = 0.5$.

Next we consider the case that $V$ is normally distributed and we use Hermite polynomials \cite{XiuBook} and the Gauss-Hermite quadrature points \cite{Gottlieb}. Let $V_{a} = \alpha, \; V_{b} = \beta$ and $V \in \left[\alpha, \; \beta \right] $, $\xi \in \left[-1, \; 1 \right]$, $\gamma \in \left(-\infty, \; \infty \right) $. The linear transformation between $V$ and $\xi$ is given by
$$ V(\xi) = \left(\frac{\beta - \alpha}{2} \right)\xi + \frac{1}{2}(\alpha + \beta) $$
and the transformation between $\xi$ and $\gamma$ is given by \cite{Chen_Gottlieb}
\begin{eqnarray}
\gamma &=& \frac{\xi}{1-\xi^2}, \qquad \xi \ne 0, \nonumber \\
&=& 0, \qquad \qquad \xi = 0. \nonumber
\end{eqnarray}
Or we have,
\begin{eqnarray}
\xi &=& \frac{-1+\sqrt{1+4\gamma^2}}{2\gamma}, \qquad \gamma \ne 0, \nonumber \\
&=& 0, \qquad \qquad \qquad \gamma = 0. \nonumber
\end{eqnarray}
Thus we have,
\begin{equation}
V(\gamma) = \left(\frac{\beta - \alpha}{2} \right)\left[\frac{-1+ \sqrt{1 + 4\gamma^2}}{2\gamma} \right]+\frac{1}{2}(\alpha + \beta), \nonumber  
\label{trans1}
\end{equation}
where $\gamma $ has the normal distribution with mean $0$ and the standard deviation (SD) $0.1$.  
For the simulation we consider $\epsilon = 0.3$ and $V \sim N\left[0,\;0.1\right]$. The figure in the bottom panel of Figure \ref{fig:soliton4} shows the mean solution at the final time obtained by the Hermite chaos. Although the mean solutions obtained from the Legendre and Hermite chaos are different, we observe that the location of $L^{'}$ is same for both cases.

Using a series of those simulations above for different values of $\epsilon$ where $\epsilon \in \left[0.05, \; 4.5\right]$,  we determine the critical velocities with respect to different  $\epsilon$. 
The results are plotted in semi-logarithmic scale  in Figure \ref{fig:soliton5}. It is observed that for the small values of $\epsilon$ where $\epsilon < 0.1$, the curve is very stiff and the slope changes sharply around $\epsilon =0.1$. From $\epsilon > 0.1$, the curve increases steadily. 
The \textit{``trapped"} and the \textit{``untrapped"} regions are clearly shown in the figure. The $V - \epsilon$ graph is the boundary of those two regions.

{{
\begin{table}
\label{table1}
\caption{Convergence of $V_{c}$ with $N$ for the Legendre Chaos. $\epsilon = 0.3,\;1.0$ and $4.5$.}
\begin{center}
\begin{tabular}
{cccc}\hline\em $N$ &\em  &  $V_{c}  \times 10^{3}$ & \\\hline &  $\epsilon = 0.3 $ & $\epsilon = 1.0$ & $\epsilon = 4.5$\\\hline\hline $2$ & $1.658675134594813$ &$23.83399810435849$&$233.3998104358486$ \\\hline $4$ & $1.786459011090171$ &$24.00997282851472$& $236.0392987896415$  \\\hline $8$ & $1.788091882999389$ &$24.01306953886891$& $236.7359886666223$ \\\hline $12$ & $1.788112748469627$ & $24.01312403314395$ &$236.9198791482719$ \\\hline $16$ & $1.788113015096694$ &$24.01312499210545$&  $236.9684168267407$ \\\hline  $20$ & $1.788113018503758$&$24.01312500898075$& $ 236.9812282904446$\\\hline $24$ & $1.788113018547295$ &$24.01312500927771$ & $236.9846098613687$  \\\hline  
\end{tabular}
\end{center}
\end{table}
}}

\subsection{Convergence analysis}
\begin{figure}
	\centering
		\includegraphics[width=0.47\textwidth]{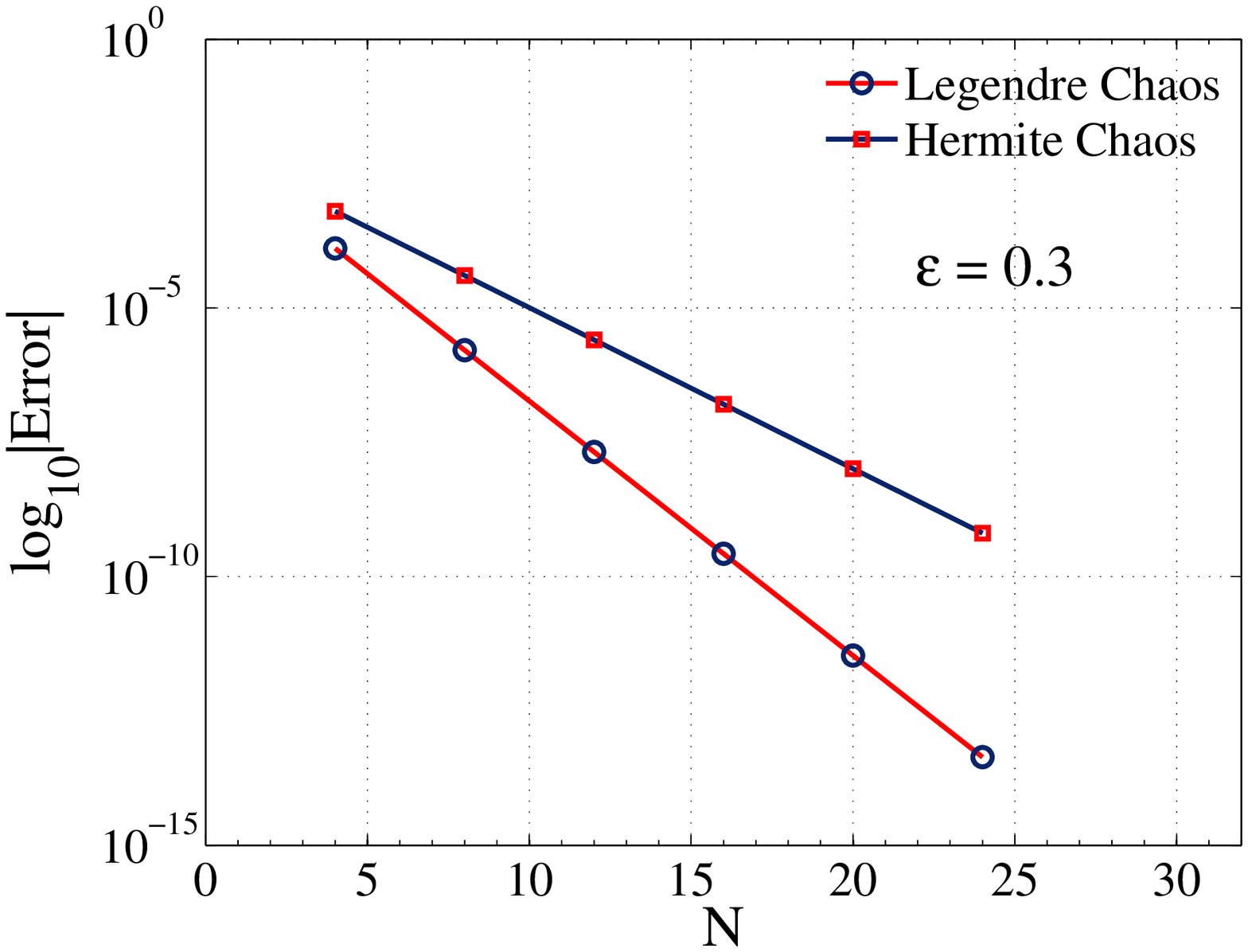}
		\includegraphics[width=0.47\textwidth]{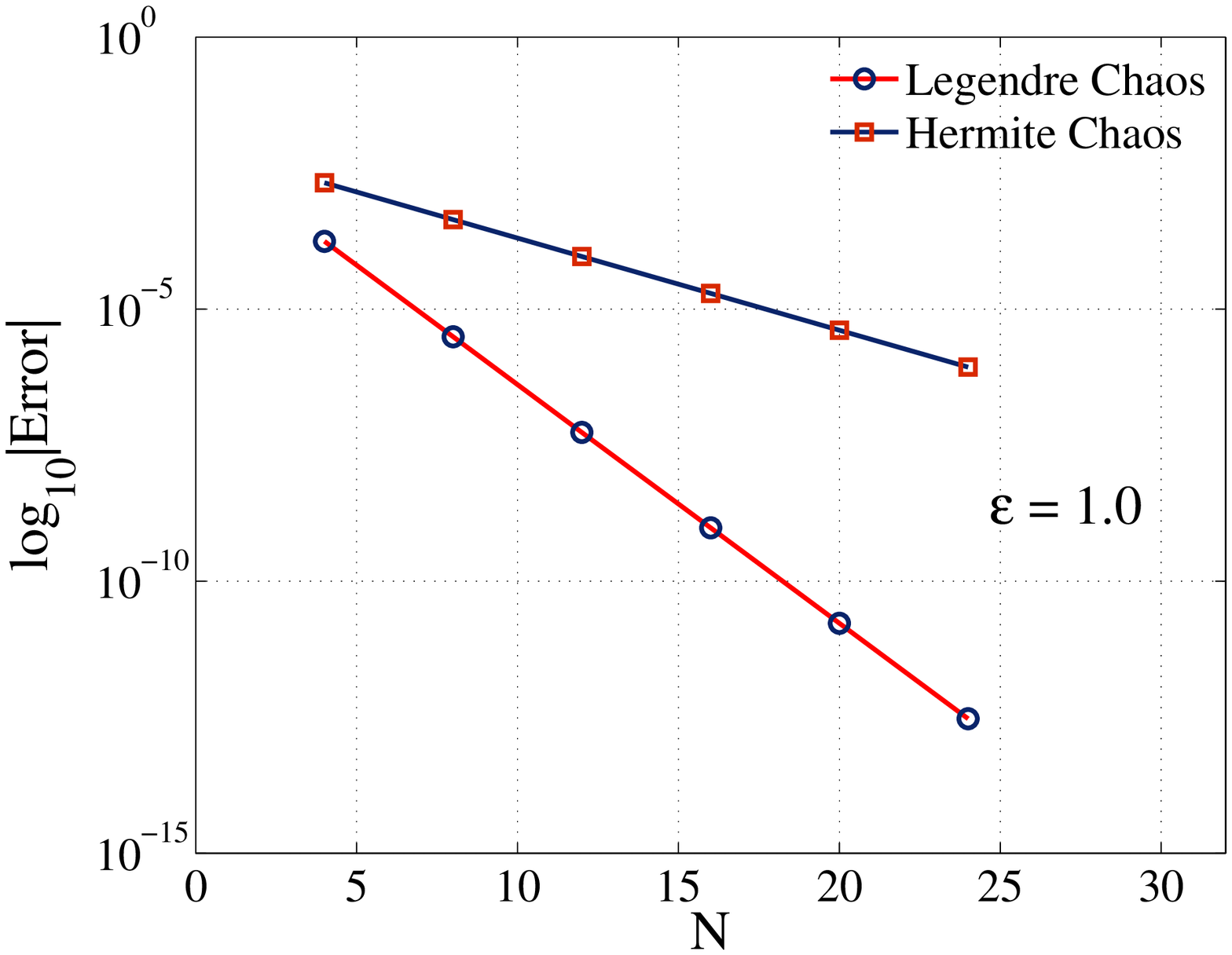}
		\includegraphics[width=0.47\textwidth]{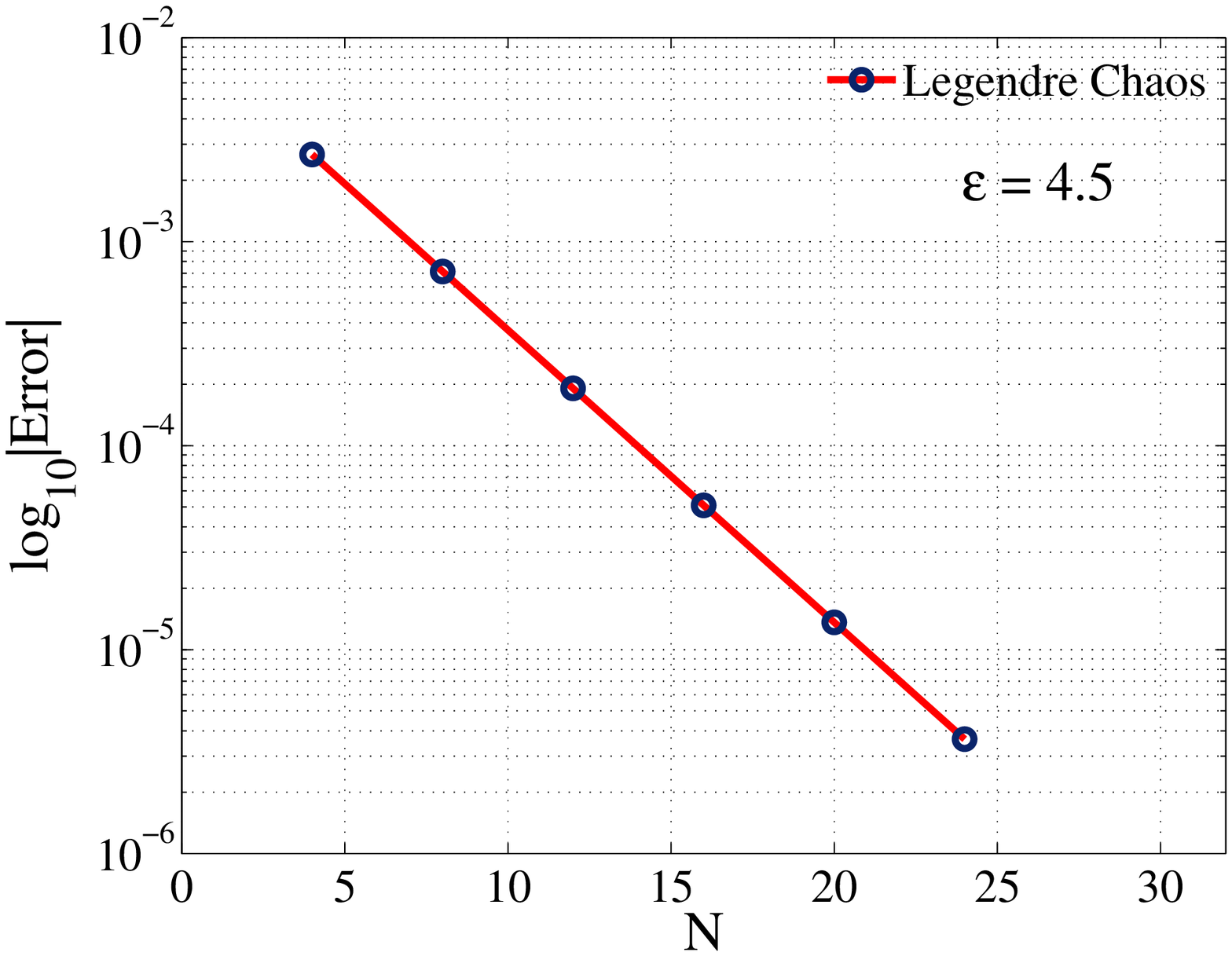}
		\caption{Spectral convergence of the critical velocities for $\epsilon = 0.3, 0.1$ and $4.5$. Graph shows the spectral convergence for both the Legendre and Hermite chaos. Note that the Legendre chaos shows faster convergence than Hermite chaos.}
	 \label{fig:soliton6}
\end{figure}

We define the error of the critical velocities  by 
$$
\mathrm{Error}^{\epsilon}(N) = \left|V_{c}^{\epsilon}(N) - V_{c}^{\epsilon}(N-1)\right|,
$$
where $N$ is the number of collocation points. Figure \ref{fig:soliton6} shows the convergence of errors obtained  by the Legendre and  Hermite chaos. We do  the convergence analysis for various values of $\epsilon$. We choose $\epsilon = 0.3$ (small) , $\epsilon = 1.0$ (moderate) and $\epsilon = 4.5$ (high).
For the Legendre chaos, the critical velocities for different $N$ are presented in Table $1$. For $\epsilon = 0.3$ and $\epsilon = 1.0$, we calculate the errors for both the Legendre and Hermite chaos and for $\epsilon = 4.5$ we use the Legendre chaos. For Hermite chaos, we expect to have the similar results. The graphs are plotted in semi-logarithmic scale. Figure \ref{fig:soliton6} shows all the graphs are a straight line, which confirms  spectral convergence but the convergence rates are different for different cases. For $\epsilon = 0.3$ and $\epsilon = 1.0$, Hermite chaos exhibits  slower convergence rate than the Legendre chaos. Also if we compare the graphs for the Legendre chaos for different cases, it is found that the convergence rate decreases with the increases of the value of $\epsilon$. That is, the smaller is the value of $\epsilon$, the faster convergence is obtained. One of the possible reasons is because of the radiation effect. As $\epsilon$ increases, the radiation effect becomes visible and it  makes difficult to locate the position of $L^{'}$ accurately. According to our numerical results, our main result is stated by the following: 
\textit{The numerical scheme stated in Section $4$ to find the critical velocity $V_{c}$ has the spectral convergence and the rate of convergence decreases with increase of the value of $\epsilon$.}

\section{Conclusion}
In this paper we studied the NLSE with the singular potential. 
We proposed an efficient method of determining the critical soliton velocities, $V_{c}$, by using the gPC collocation method. We studied the wide range of $\epsilon$, i.e. $\epsilon \in \left[0.05,\; 4.5\right]$. For $\epsilon < 0.05$ the numerical simulations demand a huge computational time due to the very small soliton velocity $\left(V \sim 10^{-10}\right)$.
We studied the convergence analysis to prove the merit of our proposed numerical scheme. We found the spectral convergence in all cases. The main development of this paper is the use of the gPC collocation method to determine the critical velocity of the soliton for given $\epsilon$ with the desired level of accuracy. We obtained $V_c$ accurately with a small number of simulations. In our future work, we will further study the case that $\epsilon \ll 0.05$. Also for the high values of $\epsilon$, where radiation effect is prominent and the convergence of the proposed method becomes slower due to the radiation effect, an efficient numerical method dealing with this effect will be investigated. 
\vskip .1in
{\bf Acknowledgement:}
\noindent
The first author is grateful to Gino Biondini for developing and implementing high-order SSFM.

\end{document}